\documentclass[11pt]{article}
\usepackage{fullpage}

\usepackage[utf8]{inputenc} % allow utf-8 input
\usepackage[T1]{fontenc}    % use 8-bit T1 fonts
\usepackage{hyperref}       % hyperlinks
\usepackage{url}            % simple URL typesetting
\usepackage{booktabs}       % professional-quality tables
\usepackage{amsfonts}       % blackboard math symbols
\usepackage{nicefrac}       % compact symbols for 1/2, etc.
\usepackage{microtype}      % microtypography
\usepackage{tcolorbox}
\usepackage{diagbox}

\usepackage{enumerate}
\usepackage[shortlabels]{enumitem}
\usepackage{algorithmic}
\usepackage{algorithm}
\usepackage{subfigure}
\usepackage{graphicx} % more modern
\usepackage{caption}
\usepackage{amsmath}
\usepackage{amsthm}
\usepackage{amssymb}
\usepackage{tikz}
\usepackage{tablefootnote}
\usepackage{multirow}
\usepackage{enumerate}
\usepackage{color}
\usepackage{xcolor}
\usepackage{natbib}

\usetikzlibrary{arrows}

\allowdisplaybreaks[4]

\usepackage{mathrsfs}

\usepackage{hyperref}
\usepackage{bm,todonotes}

\allowdisplaybreaks

\hypersetup{
colorlinks=true,
filecolor=blue,
citecolor = blue,
urlcolor=cyan,
}

\newcommand{\cO}{\mathcal{O}}

\newcommand{\RR}{\mathbb{R}}

\newcommand{\xb}{\mathbf{x}}
\newcommand{\yb}{\mathbf{y}}
\newcommand{\zb}{\mathbf{z}}
\newcommand{\TB}{\mathbb{T}}

\newcommand{\bs}{\mathbf{s}}

\newcommand{\bPi}{\mathbf{\Pi}}
\newcommand{\bbx}{\overline{x}}
\newcommand{\bbs}{\overline{s}}
\newcommand{\bby}{\overline{y}}
\newcommand{\bbz}{\overline{z}}

\newcommand{\xbti}{\xb_t^{(i)}}

\newcommand{\zbti}{\zb_t^{(i)}}

\newcommand{\argmin}{\mathop{\mathrm{argmin}}}

\newcommand{\norm}[1]{\left\|#1\right\|}
\newcommand{\dotprod}[1]{\left\langle #1\right\rangle}

\newtheorem{lemma}{Lemma}
\newtheorem{theorem}{Theorem}
\newtheorem{proposition}{Proposition}

\newtheorem{corollary}[theorem]{Corollary}
\newtheorem{remark}{Remark}
\ifx\assumption\undefined
\newtheorem{assumption}{Assumption}

\newcommand{\fm}{\mathrm{FastMix}}
\newcommand{\prox}{\mathrm{prox}}
\newcommand{\hzb}{\widehat{\zb}}

\begin{document}
\title{
Optimal Decentralized Composite Optimization for Convex Functions
}

\author{
	Haishan Ye\thanks{Center for Intelligent Decision-Making and Machine Learning, School of Management, Xi’an Jiaotong University; email: \texttt{yehaishan@xjtu.edu.cn}. }  \\
	\and
	Xiangyu Chang \thanks{Center for Intelligent Decision-Making and Machine Learning, School of Management, Xi’an Jiaotong University; email: \texttt{xiangyuchang@xjtu.edu.cn}. } 
}
\date{\today}

\maketitle

\begin{abstract}
In this paper, we focus on the decentralized composite optimization for  convex functions.
Because of advantages such as  robust to the network and no communication bottle-neck in the central server, the decentralized optimization has attracted much research attention in signal processing, control, and optimization
communities. 
Many optimal algorithms have been proposed for the objective function is  smooth and (strongly)-convex in the past years.
However, it is still an open question whether one can design an optimal algorithm when there is a non-smooth regularization term.
In this paper, we fill the gap between  smooth decentralized optimization and decentralized composite optimization and propose the first algorithm which can achieve both the optimal computation and communication  complexities.
Our experiments also validate the effectiveness and efficiency of our algorithm both in computation and communication.
\end{abstract}

\section{Introduction}
In this paper, we consider the decentralized composite optimization problem. 
In this setting, the objective function is composed of $m$ local functions
$f_i(x)$ that are located on $m$ different agents.  
These agents form a connected and undirected \emph{static} network 
and each of them only accesses its local function and communicates with its neighbors. 
All of the agents aim to cooperatively solve the following composite convex problem:
\begin{equation}
	F(x) = f(x) + g(x), \mbox{ with } f(x) = \frac{1}{m} \sum_{i=1}^{m}f_i(x), \label{eq:prob}
\end{equation}
where each $f_i(x)$ is $L$-smooth and convex, and the regularization term $g(x)$ is $\mu$-strongly convex but may be non-differentiable with $0\leq\mu\leq L$.
A typical example of $g(x)$ is the elastic-norm \citep{zou2005regularization}. 
Problems of the form~\eqref{eq:prob} arise in different settings including, among others, in machine learning \citep{wu2009genome}, model fitting \citep{boyd2011distributed}, and economic dispatch problems in power systems \citep{dominguez2012decentralized}.

Because of the advantages  such as  robust to the network and no communication bottleneck in the central server \citep{lian2017can},    decentralized optimization has attracted much research attention. 
Decentralized methods for
solving Problem~\eqref{eq:prob} have also been widely studied in signal processing, control, and optimization
communities \citep{sayed2014adaptation,shi2015extra,AlghunaimYS19,li2020revisiting}.

Due to the wide studies, lots of decentralized optimization methods have been proposed in the past years. 
For the setting that $g(x)$ is also differentiable, several algorithms can achieve the linear convergence rate \citep{shi2015extra,MokhtariSLR16,YuanYZS19,li2019decentralized,qu2017harnessing,qu2019accelerated,di2016next,li2021accelerated,kovalev2020optimal,song2023optimal,Ye2023,scaman2017optimal}.  
All these algorithms base on the primal-dual method or use the gradient-tracking method.
Among these algorithms, \citet{shi2015extra} propose the first decentralized optimization algorithm that can achieve the linear convergence rate without increasing communication cost iteratively which depends on the target precision $\epsilon$. 
\citet{scaman2017optimal} propose the first communication optimal algorithm, that is, its communication complexity is $\cO\left(\sqrt{\frac{L}{\mu(1 - \lambda_2(W))}}\log \frac{1}{\epsilon}\right)$ to find an $\epsilon$-precision solution, where $W$ is the weight matrix related to the network (refer to Assumption~\ref{ass:w}) and $\lambda_2(W)$ is the second largest eigenvalue of $W$.
\citet{kovalev2020optimal} propose the first algorithm which can achieve both computation and communication complexities, that is, its computation complexity is $\cO\left(\sqrt{\frac{L}{\mu}}\log \frac{1}{\epsilon}\right)$. 
However, the method of \citet{kovalev2020optimal} requires multi-consensus, that is, an inner communication loop is required, to achieve the optimal complexities. 
Recently, \citet{song2023optimal} propose an optimal decentralized algorithm in which each agent only communicates with its neighbors once for each iteration.

For the case that $g(x)$ is non-differentiable, there are also many works. 
Many gradient tracking based algorithms have been extended to decentralized composite optimization problems with a non-smooth regularization term such as \texttt{PG-EXTRA} \citep{ShiLWY15} and \texttt{NIDS} \citep{li2019decentralized}.
However, due to the non-smooth term, the theoretical analysis of these algorithms can only guarantee  sub-linear convergence rates.
Recently, \citet{AlghunaimYS19} proposed a primal-dual algorithm which can achieve a linear convergence rate.
At almost the same time, \cite{sun2022distributed} propose a gradient
tracking based method called \texttt{SONATA}, and established a linear convergence rate. 
However,  \texttt{SONATA} requires  multi-consensus steps and its inner communication number of multi-consensus step depends on $\log\frac{L}{\mu}$.
After the work of  \citet{AlghunaimYS19}, a unified framework is proposed to analyze a large
group of primal-dual and gradient tracking based algorithms, and show that these algorithms can also achieve linear convergence rates with nonsmooth regularization term such as \texttt{EXTRA} (\texttt{PG-EXTRA}) \citep{shi2015extra}, and \texttt{NIDS} \citep{li2019decentralized} in the work of \citep{alghunaim2020decentralized,xu2021distributed}.
Despite intensive studies in the literature, the convergence rates of
these previous algorithms do not match the optimal convergence rate. 
Moreover, the communication complexities achieved by algorithms
analyzed in the framework of \citet{xu2021distributed} and \citet{alghunaim2020decentralized} are sub-optimal.
Though several optimal algorithms have been proposed when the regularization term $g(x)$ is smooth, it still remains open whether one can design an optimal algorithm for the the decentralized composite optimization.
To the best of our knowledge, \texttt{ProxMudag} is the best of proposed algorithms for the decentralized composite optimization for the strongly convex function \citep{Ye2023}.
It can achieve the optimal computation complexity but only achieve a near optimal communication complexity $\cO\left(\sqrt{\frac{L}{\mu(1 - \lambda_2(W))}}\log\frac{L}{\mu}\log \frac{1}{\epsilon}\right)$ with an extra $\log\frac{L}{\mu}$ factor compared to the optimal communication complexity.
On the other hand, all works mentioned above focus on  strongly convex functions. 
When $f_i(x)$ and $g(x)$ are generally convex, much less works have studied this case.
Thus, it also remain open whether one can design an optimal algorithm for  the decentralized composite optimization for general convex functions.

In this paper, we try to fill the gap between the smooth decentralized optimization and decentralized composite optimization. 
We propose the first algorithm that can achieve both the optimal computation and communication  complexities when  $f_i(x)$ is smooth and the regularization term $g(x)$ is (strongly)-convex. 
Furthermore, our method can be easily extended for the problems that $f_i(x)$ is $L$-smooth and $\mu$-strongly convex while $g(x)$ is convex. 
Our algorithm can still achieve both the optimal computation and communication  complexities (refer to Section~\ref{subsec:extension}). 
Our algorithm mainly consists of Nesterov's acceleration, gradient tracking, and multi-consensus.
Thus it is easy to implement and its parameters are easy to tune.
Our experiments also validate the effectiveness and efficiency of our algorithm both in computation and communication.

\section{Notation and Assumptions}

Let $\xb$, $\yb$, $\zb$ and $\bs$ are $m\times d$ matrices whose $i$-th rows $\xb^{(i)}$, $\yb^{(i)}$, $\zb^{(i)}$ and $\bs^{(i)}$ are their  local copies for the $i$-th agent, respectively. 
Accordingly, we define the averaging variables 
\begin{equation}\label{eq:bbx}
	\begin{aligned}
		\bbx :=& \frac{1}{m}\sum_{i=1}^m \xb^{(i)} = \frac{1}{m} \mathbf{1}^\top \xb \in\RR^{1\times d}, \quad \bby:= \frac{1}{m} \sum_{i=1}^{m} \yb^{(i)}  \\
		\bbz :=& \frac{1}{m} \sum_{i=1}^{m}\zb^{(i)} \quad \bbs := \frac{1}{m}\mathbf{1}^\top \bs \in\RR^{1\times d},
	\end{aligned}
\end{equation}
where $\mathbf{1}$ denotes the vector with all entries equal to $1$ of dimension $m$. 
Now we introduce the projection matrix
\begin{equation}\label{eq:Pi_def}
	\bPi = \mathbf{I}_m - \frac{\mathbf{1}\mathbf{1}^\top}{m}.
\end{equation}
Using the projection matrix $ \bPi $, we can represent that
\begin{equation*}
	\norm{\xb - \mathbf{1}\bbx} = \norm{ \xb - \frac{\mathbf{1}\mathbf{1}^\top}{m} \xb  } = \norm{\bPi \xb}.
\end{equation*}

Now, we introduce an aggregate objective function $\widetilde{f}(\xb)$ defined as follows,
\begin{equation}
	{f}(\xb) := \sum_{i=1}^m f_i(\xb^{(i)}).
\end{equation}
Accordingly, its aggregate gradient
\begin{equation*}
	\nabla {f}(\xb) := [\nabla f_1(\xb^{(1)}); \nabla f_2(\xb^{(2)});\dots; \nabla f_m(\xb^{(m)})] \in\RR^{m\times d}.
\end{equation*}
It holds that the $i$-th row of $\nabla \widetilde{f}(\xb)$ equals to $\nabla f_i(\xb^{(i)})$.

Throughout this paper, we use $ \norm{\cdot} $ to denote the ``Frobenius'' norm. 
That is, for a matrix $ \xb \in \RR^{m\times d} $, it  holds that
\begin{align*}
	\norm{\xb}^2 = \sum_{i=1, j=1}^{m, d} \left(\xb^{(i, j)}\right)^2,
\end{align*}
where $\xb^{(i,j)}$ denote the entry of $\xb$ in the $i$-th row and $j$-th column.
Furthermore, we use $ \norm{\xb}_2 $ to denote the spectral norm which is the largest singular value of $ \xb $.
For vectors $ x, y \in \RR^d $, we use $ \dotprod{x, y} $ to denote the standard inner product of $ x $ and $ y $.

Next, we also introduce the aggregated proximal operator with respect to $g(\cdot)$ as for $\xb\in\RR^{m\times d}$,
\begin{align*}
	\prox_{\gamma g}(\xb) = \argmin_{\mathbf{w}\in\RR^{m\times d}} \left( \gamma g(\mathbf{w}^{(i)}) + \frac{1}{2} \norm{ \mathbf{w}^{(i)} - \xb^{(i)} }^2 \right).
\end{align*} 
Thus, the $i$-th row of $\prox_{\gamma g}(\xb)$ equals to $\prox_{\gamma g}(\xb^{(i)})$.

Now we introduce several assumptions that will be used throughout this paper.
We focus on the smooth and strongly-convex functions. 
That is, the function $f_i(x)$ satisfies the following assumption.
\begin{assumption}\label{ass:f}
	Each $f_i:\RR^d\to \RR$ is convex  and $L$-smooth, i.e., for any $x, y \in \RR^d$, it holds that
	\begin{align*}
		f_i(y) \ge& f_i(x) + \dotprod{\nabla f_i(x), y -x } ,\\
		f_i(y) \le& f_i(x) + \dotprod{\nabla f_i(x), y -x } + \frac{L}{2}\norm{x - y}^2.
	\end{align*}
\end{assumption}

\begin{assumption}\label{ass:g}
	The function $g: \RR^d \to \RR$ is $\mu$-strongly convex, i.e., for any $x, y \in\RR^d$, it holds that
	\begin{equation*}
		g(y) \ge g(x) + \dotprod{\partial g(x), y - x} + \frac{\mu}{2}\norm{y - x}^2, \quad\mbox{where}\quad \partial g(x) \mbox{ is the subgradient at } x.
	\end{equation*}
\end{assumption}

Using the $f(\cdot)$, $f_i(\cdot)$, and aggregated variable $\xb$, we define the extended Bregman distance between $y\in\RR^d$ and $\xb\in\RR^{m\times d}$ as follows:
\begin{equation}\label{eq:df}
	D_f(y, \xb) \triangleq \frac{1}{m} \sum_{i=1}^{m} \left( f_i(y) - f_i(\xb^{(i)}) - \dotprod{\nabla f_i(\xb^{(i)}), y - \xb^{(i)}} \right).
\end{equation}
Because of the convexity of $f_i(\cdot)$, $D_f(y, \xb)$ is non-negative for any $y\in\RR^d$ and $\xb\in\RR^{m\times d}$.

For the topology of the network, we let $W$ be the weight matrix associated with the network, indicating how agents are connected to each other.
The weight matrix $W$ satisfies the following assumption.
\begin{assumption}\label{ass:w}
	The weight matrix $W$ is symmetric positive semi-definite with $W_{i,j} \neq 0$ if and if only agents $i$ and $j$ are connected or $i=j$. It also satisfies that ${\mathbf{0}}\preceq W\preceq I$, $W{\mathbf{1}}  = {\mathbf{1}}$, ${\rm null}(I - W) = {\rm span}(\mathbf{1})$.
\end{assumption}
The weight matrix satisfying the above assumption has an important property that $W^\infty = \frac{1}{m}\mathbf{1}\mathbf{1}^\top$ \citep{xiao2004fast}.
Thus, one can achieve the effect of averaging local $\xb^{(i)}$ on different
agents by using $W\xb$ for iterations.
Instead of directly multiplying $W$, \citet{liu2011accelerated} proposed a more efficient way to achieve averaging described in Algorithm~\ref{alg:mix}, which has the following important proposition.
\begin{proposition}[\cite{Ye2023}]
	\label{lem:mix_eq}
	Let $\xb^K$ be the output of Algorithm~\ref{alg:mix} with $\eta_w = 1/(1+\sqrt{1-\lambda_2^2(W)})$ and we denote~$\bar{x} = \frac{1}{m}\mathbf{1}^\top\xb^1$.
	Then it holds that 
	\[
	\bar{x} = \frac{1}{m}\mathbf{1}^\top\xb^K \quad\mbox{and}\quad 
	\norm{\xb^K - \mathbf{1}\bar x} 
	\leq \sqrt{14} \left(1 - \left(1-\frac{1}{\sqrt{2}}\,\right) \sqrt{1-\lambda_2(W)}\right)^K \norm{\xb^1 - \mathbf{1}\bar x},
	\]
	where $\lambda_2(W)$ is the second largest eigenvalue of $W$. 
\end{proposition}

\section{Optimal Decentralized Accelerated Proximal Gradient Descent}

In this section, we first give the algorithm description.
Then, we give the main results of this paper.
Finally, we extend our algorithm to solve the problem that $f_i(x)$ is $L$-smooth and $\mu$-strongly convex while $g(x)$ is convex but may be non-differentiable.
\begin{algorithm}[tb]
	\caption{Optimal Decentralized Accelerated Proximal Gradient Descent (\texttt{ODAPG})}
	\label{alg:alg_name}
	\small
	\begin{algorithmic}
		\STATE {\bfseries{Input}}: $x_1$, mixing matrix $W$,  step size sequence $\{\gamma_t\}$ and sequence $\{\tau_t\}$.
		\STATE {\bfseries{Initialization}:} Set $\zb_1 = \yb_1 = \xb_1 = \mathbf{1}x_1$,  $\bs_1^{(i)} = \nabla f_i(\xb_1^{(i)})$, in parallel for $i \in [m]$.
		\FOR {$t = 1,\dots, T$}
		\STATE Compute  $ \xb_{t+1} = \tau_t \zb_t + (1 - \tau_t) \yb_t $. 
		\STATE Compute the local gradients $\nabla f_i(\xb_{t+1}^{(i)})$ in parallel for $i\in [m]$ to form the gradient $\nabla {f}(\xb_{t+1})$. \label{step:grad}
		\STATE $\bs_{t+1} = \fm( \bs_t + \nabla  {f}(\xb_{t+1}) - \nabla {f}(\xb_t), K )$ 
		\STATE $ \zb_{t+1} = \fm\left(\prox_{\gamma g} \left(\zb_t - \gamma_t \bs_{t+1}\right),K\right) $ \label{step:prox}
		\STATE $\yb_{t+1} = \fm\left(\tau_t \zb_{t+1} + (1 - \tau_t) \yb_t, K\right)$
		\ENDFOR
	\end{algorithmic}
\end{algorithm}

\begin{algorithm}[tb]
	\caption{FastMix}
	\label{alg:mix}
	\begin{small}
		\begin{algorithmic}[1]
			\STATE {\bf Input:} $\xb^{1}$, $K$, $W$ and $\eta_w$ \\[0.15cm]
			\STATE $\xb^{0} = \xb^{1}$ \\[0.15cm]
			\STATE\textbf{for} $k=1,\dots, K$ \textbf{do} \\[0.15cm]
			\STATE\quad $\xb^{k+1} = (1+\eta_w)W\xb^k - \eta_w \xb^{k-1}$ \\[0.15cm]  
			\STATE\textbf{end for} \\[0.15cm]
			\STATE {\bf Output:} $\xb^K$
		\end{algorithmic}
	\end{small}
\end{algorithm}
\subsection{Algorithm Description}
Our algorithm combines the Nesterov' accelerated gradient descent and gradient tracking. The main algorithmic procedure of our algorithm is listed as follows:
\begin{align}
	\xb_{t+1} =& \tau_t \zb_t + (1 - \tau_t) \yb_t \label{eq:xzy} \\
	\bs_{t+1} =& \fm\big( \bs_t + \nabla \widetilde{f}(\xb_{t+1}) - \nabla \widetilde{f}(\xb_t), K \big) \label{eq:ss}\\
	\hzb_{t+1} =& \prox_{\gamma g}(\zb_t - \gamma_t \bs_{t+1}) \label{eq:hz} \\
	\zb_{t+1} =&\fm \left( \hzb_{t+1}, K\right), \label{eq:prox_step}\\
	\yb_{t+1} =& \fm \left(\tau_t \zb_{t+1} + (1 - \tau_t) \yb_t, K\right), \label{eq:yzy1}
\end{align}
where $K$ is the iteration number for the ``FastMix'' operator. 
We can observe that our algorithm almost shares the same procedure with the one of accelerated proximal gradient descent except the gradient $\nabla \widetilde{f}(\xb_{t+1})$ replaced with the gradient tracking variable $\bs_{t+1}$ and the extra ``FastMix'' operator.

The ``accelerated proximal gradient descent'', ``gradient tracking'', and ``FastMix'' are the three important elements to the success of our algorithm to achieve the optimal communication and computation  complexities. 
The ``accelerated proximal gradient descent'' provides the potential  that our algorithm can achieve a fast convergence rate.
The ``gradient tracking'' can be regarded as a kind of variance reduction. 
Its update rule is the same to the one of an important variance reduction named SARAH \citep{nguyen2017sarah}.
Thus, the ``gradient tracking'' technique is very useful to reduce the communication complexity and has been widely used in decentralized optimization.
To further reduce the communication complexity, the ``FastMix'' is the key.  
Without the ``FastMix'', the communication complexity depends on $(1 - \lambda_2(W))^{-1}$. 
In contrast, the ``FastMix'' helps to achieve a communication complexity depending on $(1 - \lambda_2(W))^{-1/2}$.

\subsection{Main Results}

In this work, we focus on the synchronized setting in which the computation complexity depends
on the number of gradient calls and the communication complexity depends on the rounds of local
communication. We give the detailed upper complexity bounds of our algorithm when $f(x)$ is $L$-smooth and $g(x)$ is $\mu$-strongly convex  in the following
theorem.

\begin{theorem}\label{thm:main}
	Let the objective function $F(x)$ be of the form~\eqref{eq:prob} and Assumption~\ref{ass:f}-\ref{ass:g} hold with $\mu>0$. 	
	Letting us set a constant step size $\gamma = \frac{1}{20\sqrt{L\mu}}$, $\tau = \mu \gamma$, and $K = \frac{15}{\sqrt{1 - \lambda_2(W)}}$ in Algorithm~\ref{alg:alg_name}, letting $x^*$ denote the minimum of $F(x)$, then it holds that
	\begin{equation}\label{eq:main}
		\begin{aligned}
			\norm{\zb_T -\mathbf{1}x^*}^2
			\leq& 
			\left(1 - \frac{1}{40}\cdot \sqrt{\frac{\mu}{L}}\right)^T  \cdot \Bigg( \frac{2m}{\mu} \big(F(\bby_1) - F(x^*)\big) +  \norm{\zb_1 - \mathbf{1}x^*}^2 \\
			&+ \frac{ 6\cdot20^2 L}{\mu}\cdot \left( \norm{\Pi\xb_1}^2 + c_1 \norm{\Pi\yb_1}^2 + c_2 \norm{\Pi \zb_1}^2 + c_3\gamma^2 \norm{\Pi\bs_1}^2 \right)  \Bigg),
		\end{aligned}
	\end{equation}
	with $c_1 = 64$, $c_2 = \frac{9\tau^2}{32}$, and $c_3 = \frac{5\tau^2}{14}$. 
\end{theorem}

Above theorem shows that our algorithm converges with a rate $1 - \cO\left(\sqrt{\frac{\mu}{L}}\right)$ which is the same to the one of accelerated proximal gradient descent. 
Furthermore, since it holds that $ \norm{\zbti - x^*}^2 \leq \norm{\zb_t - \mathbf{1}x^*}^2$, Eq.~\eqref{eq:main} shows that $\zbti$ in $i$-th agent will converge to the optimum.

Next, by the convergence rate of our algorithm shown in Theorem~\ref{thm:main}, we will give the computation cost and communication complexities of our algorithm in the following corollary.
\begin{corollary}\label{cor:main}
	Let the objective function satisfy the property described in Theorem~\ref{thm:main}. 
	To find an $\epsilon$-suboptimal solution, the computation and communication complexities of Algorithm~\ref{alg:alg_name} for each agent are 
	\begin{align}
		T = \cO\left( \sqrt{\frac{L}{\mu}} \log\frac{1}{\epsilon}  \right),\quad\mbox{and}\quad C = \cO\left( \sqrt{\frac{L}{\mu (1 - \lambda_2(W))}} \log \frac{1}{\epsilon} \right). \label{eq:complexity}
	\end{align}
\end{corollary}
\begin{proof}
	By Theorem~\ref{thm:main}, we can obtain that
	\begin{align*}
		\norm{\zb_T - \mathbf{1}x^*}^2
		\leq&
		\left(1 - \frac{1}{40}\cdot \sqrt{\frac{\mu}{L}}\right)^T  \cdot \Bigg( \frac{2m}{\mu} \big(F(\bby_1) - F(x^*)\big) +  \norm{\zb_1 - \mathbf{1}x^*}^2 \\
		&+ \frac{ 6\cdot20^2 L}{\mu}\cdot \left( \norm{\Pi\xb_1}^2 + c_1 \norm{\Pi\yb_1}^2 + c_2 \norm{\Pi \zb_1}^2 + c_3\gamma^2 \norm{\Pi\bs_1}^2 \right)  \Bigg)\\
		\leq&
		\exp\left(- \frac{1}{40}\cdot \sqrt{\frac{\mu}{L}} T\right) \cdot \Bigg( \frac{2m}{\mu} \big(F(\bby_1) - F(x^*)\big) +  \norm{\zb_1 - \mathbf{1}x^*}^2 \\
		&+ \frac{ 6\cdot20^2 L}{\mu}\cdot \left( \norm{\Pi\xb_1}^2 + c_1 \norm{\Pi\yb_1}^2 + c_2 \norm{\Pi \zb_1}^2 + c_3\gamma^2 \norm{\Pi\bs_1}^2 \right)  \Bigg).
	\end{align*}
	Thus, to achieve that $\norm{\zbti - x^*}^2 \leq \norm{\zb_T - \mathbf{1}x^*}^2 \leq \epsilon $, $T$ is required to be
	\begin{align*}
		T =& 40\sqrt{\frac{L}{\mu}} \Bigg(\log \frac{1}{\epsilon} +  \log \Big( \frac{2m}{\mu} \big(F(\bby_1) - F(x^*)\big) +  \norm{\zb_1 - \mathbf{1}x^*}^2 \\
		&+  \frac{ 6\cdot20^2 L}{\mu}\cdot \left( \norm{\Pi\xb_1}^2 + c_1 \norm{\Pi\yb_1}^2 + c_2 \norm{\Pi \zb_1}^2 + c_3\gamma^2 \norm{\Pi\bs_1}^2 \right)  \Big)\Bigg)\\
		=& \cO\left( \sqrt{\frac{L}{\mu}}\log \frac{1}{\epsilon} \right).
	\end{align*}
	As shown in Algorithm~\ref{alg:alg_name}, each agent computes its local gradient $\nabla f_i(\xb_{t+1}^{(i)})$ for each iteration. Thus, the computation complexity of our algorithm equals to $T$. 
	
	Our algorithm takes three ``FastMix'' steps in which each agent communicates $K$ times with its neighbors  for each iteration.    The communication complexity of our algorithm is 
	\begin{align*}
		C = 3TK = T \cdot \frac{45}{\sqrt{1 - \lambda_2(W)}} = \cO\left( \sqrt{\frac{L}{\mu (1 - \lambda_2(W))}} \log \frac{1}{\epsilon} \right). 
	\end{align*}
\end{proof}

In the following theorem, we will give the detailed upper complexity bounds of our algorithm when $f(x)$ is $L$-smooth and $g(x)$ is general convex.

\begin{theorem}\label{thm:main1}
Let the objective function $F(x)$ be of the form~\eqref{eq:prob} and Assumption~\ref{ass:f}-\ref{ass:g} hold with $\mu=0$, that is, $g(x)$ is only convex. 	
Denoting $c_f = 200$ and setting step size $\gamma_t = \frac{t+4}{2Lc_f}$, $\tau_t = \frac{1}{Lc_f\gamma_t} = \frac{2}{t+4}$, and $K = \frac{15}{\sqrt{1 - \lambda_2(W)}}$ in Algorithm~\ref{alg:alg_name}, letting $x^*$ denote the minimum of $F(x)$, then it holds that
\begin{equation}
\begin{aligned}
F(\bby_T) 
\leq& 
\frac{15}{(T+3)^2} \Big(F(\bby_1) - F(x^*)\Big) + \frac{2Lc_f}{m(T+3)^2} \norm{\zb_1 - \mathbf{1}x^*}^2 \\
&
+ \frac{50}{(T+3)^2} \left(\norm{\Pi\xb_1}^2 + \frac{9\tau_1^2}{32}\norm{\Pi\zb_1}^2 + \frac{5}{14L^2c_f^2}\norm{\Pi\bs_1}^2 \right).
\end{aligned}
\end{equation}
\end{theorem}

\begin{corollary}
Let the objective function satisfy the property described in Theorem~\ref{thm:main1}. 
To find an $\epsilon$-suboptimal solution, the computation and communication complexities of Algorithm~\ref{alg:alg_name} for each agent are 
\begin{align}
	T = \cO\left( \frac{1}{\sqrt{\epsilon}}  \right),\quad\mbox{and}\quad C = \cO\left( \frac{1}{\sqrt{\epsilon (1 - \lambda_2(W))}}  \right). \label{eq:complexity1}
\end{align}
\end{corollary}
\begin{proof}
By Theorem~\ref{thm:main1}, to find an $\epsilon$-suboptimal solution, Algorithm~\ref{alg:alg_name} only takes 
\begin{equation*}
	T = \cO\left(\frac{1}{\epsilon^{1/2}}\right).
\end{equation*}

Our algorithm takes three ``FastMix'' steps in which each agent communicates $K$ times with its neighbors  for each iteration.    
The communication complexity of our algorithm is 
\begin{align*}
	C = 3TK = T \cdot \frac{45}{\sqrt{1 - \lambda_2(W)}} = \cO\left( \frac{1}{\sqrt{\epsilon (1 - \lambda_2(W))}}  \right). 
\end{align*}
\end{proof}

\begin{remark}
Eq.~\eqref{eq:complexity} and Eq.~\eqref{eq:complexity1} show that our algorithm can achieve optimal computation complexities  which match the computation complexity of accelerated proximal gradient descent both for the strongly convex and general convex functions.  
Thus, our algorithm achieves the optimal computation complexities \citep{nesterov2003introductory}.
Furthermore, the communication complexities of our algorithm  match the lower bounds of the decentralized optimization over strongly and general convex functions \citep{scaman2017optimal}.
Thus, our algorithm achieves the optimal computation and communication complexities.
\end{remark}

\subsection{Extension}
\label{subsec:extension}

Corollary~\ref{cor:main} shows that our algorithm can achieve the optimal computation and communication complexities for the case that $f_i(x)$ is convex and $L$-smooth while $g(x)$ is $\mu$-strongly convex.
Next, we will show that our algorithm can be extend to solve this kind of problems that $f_i(x)$ is $\mu$-strong convex and $L$-smooth while $g(x)$ is only convex. 
At the same time, our algorithm can also achieve the optimal computation and communication complexities.

\begin{algorithm}[tb]
	\caption{\texttt{ODAPG} for $f_i(x)$ being $L$-smooth and $\mu$-strongly convex, and $g(x)$ is convex.}
	\label{alg:alg_name1}
	\small
	\begin{algorithmic}
		\STATE {\bfseries{Input}}: $x_1$, mixing matrix $W$, initial step size $\widehat{\gamma}$, $\widehat{\tau}$.
		\STATE {\bfseries{Initialization}:} Set $\zb_1 = \yb_1 = \xb_1 = \mathbf{1}x_1$,  $\bs_1^{(i)} = \nabla f_i(\xb_1^{(i)}) - \mu x_1$, in parallel for $i \in [m]$.
		\FOR {$t = 1,\dots, T$}
		\STATE Compute  $ \xb_{t+1} = \widehat{\tau} \zb_t + (1 - \widehat{\tau}) \yb_t $. 
		\STATE Compute the local gradients $\nabla \widehat{f}_i(\xb_{t+1}^{(i)}) = \nabla f_i(\xb_{t+1}^{(i)}) - \mu \xb_{t+1}^{(i)}$ in parallel for $i\in [m]$ to form the gradient $\nabla \widehat{f}(\xb_{t+1})=\left[\nabla \widehat{f}_i(\xb_{t+1}^{(1)}); \nabla \widehat{f}_i(\xb_{t+1}^{(2)});\dots;\nabla \widehat{f}_i(\xb_{t+1}^{(m)})\right]$.
		\STATE $\bs_{t+1} = \fm( \bs_t + \nabla  \widehat{f}(\xb_{t+1}) - \nabla \widehat{f}(\xb_t), K )$
		\STATE $ \zb_{t+1} = \fm \Big(\prox_{\frac{\widehat{\gamma}}{1+\mu\widehat{\gamma}}g} \left(\frac{1}{1+\mu\widehat{\gamma}}\left(\zb_t - \widehat{\gamma} \bs_{t+1}\right)\right), K\Big) $
		\STATE $\yb_{t+1} = \fm\left(\widehat{\tau} \zb_{t+1} + (1 - \widehat{\tau}) \yb_t, K\right)$
		\ENDFOR
	\end{algorithmic}
\end{algorithm}

The main idea of our extension of our algorithm relies on the fact that  $\widehat{f}_i(x) = f_i(x) - \frac{\mu}{2}\norm{x}^2$ is $(L-\mu)$-smooth and convex, $\widehat{g}(x) = g(x) + \frac{\mu}{2}\norm{x}^2$ is $\mu$-strongly convex if $f_i(x)$ is $\mu$-strong convex and $L$-smooth while $g(x)$ is only convex. 
Then, we can use Algorithm~\ref{alg:alg_name} with a step size $\widehat{\gamma} = \frac{1}{20\sqrt{(L-\mu)\mu}}$ and $\widehat{\tau} = \mu \widehat{\gamma}$ to solve the problem represented as follow:
\begin{align*}
	F(x) = \frac{1}{m} \sum_{i=1}^{m}\widehat{f}_i(x) + \widehat{g}(x).
\end{align*}
Instead, the gradient computation in Step~\ref{step:grad} of Algorithm~\ref{alg:alg_name} becomes $$\nabla \widehat{f}(\xb_t) = \left[\nabla f_1(\xb^{(1)}); \nabla f_2(\xb^{(2)});\dots;\nabla f_m(\xb^{(m)})\right] - \mu \left[\xb_t^{(1)}; \xb_t^{(2)};\dots;\xb_t^{(m)}\right].$$
Similarly, the proximal mapping in Step~\ref{step:prox} becomes 	$\prox_{\widehat{\gamma}\widehat{g}}\left(\zb - \widehat{\gamma}\bs_{t+1}\right) $. 
At the same time, it holds that
\begin{align*}
	\prox_{\widehat{\gamma}\widehat{g}}\left(\zbti - \widehat{\gamma}\bs_{t+1}^{(i)}\right) 
	=& 
	\argmin_w \left( \widehat{g}(w) + \frac{1}{2\widehat{\gamma}} \norm{w - \left( \zbti - \widehat{\gamma}\bs_{t+1}^{(i)} \right)  }^2 \right)\\
	=&
	\argmin_w\left(g(w) + \frac{\mu}{2}\norm{w}^2 + \frac{1}{2\widehat{\gamma}}  \norm{w - \left( \zbti - \widehat{\gamma}\bs_{t+1}^{(i)} \right)  }^2\right)\\
	=& 
	\argmin_w\left(g(w) + \frac{1}{2 \left( \frac{\widehat{\gamma}}{1+\mu\widehat{\gamma}} \right)} \norm{  w - \frac{1}{\widehat{\gamma}}\frac{\widehat{\gamma}}{1+\mu\widehat{\gamma}}\left( \zbti - \widehat{\gamma}\bs_{t+1}^{(i)}\right)}^2 \right)\\
	=&
	\prox_{ \frac{\widehat{\gamma}}{1+\mu\widehat{\gamma}}g}\left(\frac{1}{1+\mu\widehat{\gamma}}\left( \zbti - \widehat{\gamma}\bs_{t+1}^{(i)}\right)\right).
\end{align*}
Thus, Step~\ref{step:prox} of Algorithm~\ref{alg:alg_name} is replaced by $\zb_{t+1} = \prox_{ \frac{\widehat{\gamma}}{1+\mu\widehat{\gamma}}g}\left(\frac{1}{1+\mu\widehat{\gamma}}\left( \zb - \widehat{\gamma}\bs_{t+1}\right)\right)$.

We give a detailed algorithm description in Algorithm~\ref{alg:alg_name1} for solving the problems that $f_i(x)$ is $L$-smooth and $\mu$-strongly convex while $g(x)$ is only convex.
We provide the convergence rate of Algorithm~\ref{alg:alg_name1} in the following theorem.

\begin{theorem}
	Let the objective function $F(x)$ be of the form~\eqref{eq:prob} and  each $f_i(x)$ be $L$-smooth and $\mu$-strongly convex with $L\ge 2\mu$. Suppose that $g(x)$ is convex but may be non-smooth.
	Let us set the step size $\widehat{\gamma} = \frac{1}{20\sqrt{(L-\mu)\mu}}$, $\widehat{\tau} = \mu \widehat{\gamma}$, and $K = \frac{11}{\sqrt{1 - \lambda_2(W)}}$ in Algorithm~\ref{alg:alg_name1}, then the output $\zb_T$ satisfies that
	\begin{equation}\label{eq:main3}
		\begin{aligned}
			&\norm{\zb_T - \mathbf{1}x^*}^2
			\leq
			\left(1 - \frac{1}{40}\cdot \sqrt{\frac{\mu}{L-\mu}}\right)^T  \cdot \Bigg( \frac{2m}{\mu} \big(F(\bby_1) - F(x^*)\big) +  \norm{\zb_1 - \mathbf{1}x^*}^2 \\
			&+ \frac{ 6\cdot20^2 (L-\mu)}{\mu}\cdot \left( \norm{\Pi\xb_1}^2 + c_1 \norm{\Pi\yb_1}^2 + c_2 \norm{\Pi \zb_1}^2 + c_3\gamma^2 \norm{\Pi\bs_1}^2 \right) \Bigg),
		\end{aligned}
	\end{equation}
	with $c_1 = 64$, $c_2 = \frac{9\widehat{\tau}^2}{32}$, and $c_3 = \frac{5\widehat{\tau}^2}{14}$.  
\end{theorem}
\begin{proof}
	Let us denote that $\widehat{g}(x) = g(x) + \frac{\mu}{2}\norm{x}^2$ and $\widehat{f}_i(x) = f_i(x) - \frac{\mu}{2}\norm{x}^2$.
	Since $f_i(x)$ is $L$-smooth and $\mu$-strongly convex, and $g(x)$ is only convex but may be not smooth, then $\widehat{g}(x)$ is $\mu$-strongly convex and $\widehat{f}_i(x)$ is $(L-\mu)$-smooth and convex which satisfy the conditions of Theorem~\ref{thm:main}. 
	Running Algorithm~\ref{alg:alg_name} with $\widehat{g}(x)$ and $\widehat{f}(x)$ with a step size $\widehat{\gamma} = \frac{1}{20\sqrt{(L-\mu)\mu}}$ and $\widehat{\tau} = \mu \widehat{\gamma}$, then Theorem~\ref{thm:main} guarantees the convergence rate in Eq.~\eqref{eq:main3}. 
\end{proof}

\begin{remark}
	By the similar analysis in the proof of Corollary~\ref{cor:main}, we can obtain Algorithm~\ref{alg:alg_name1} can achieve a computation complexity  is $\cO\left( \sqrt{\frac{L}{\mu}} \log \frac{1}{\epsilon} \right)$ and a communication complexity $\cO\left( \sqrt{\frac{L}{\mu(1-\lambda_2(W))}} \log \frac{1}{\epsilon} \right)$.
	Thus, our work provides an optimal decentralized algorithm for solving problems that $f_i(x)$ is $L$-smooth and $\mu$-strongly convex while $g(x)$ is only convex.
\end{remark}

\section{Convergence Analysis}

In this section, we will prove Theorem~\ref{thm:main}. 
First, we will provide important lemmas related to ``accelerated proximal gradient descent''.
Next, we will bound the the consensus error terms.
Finally, we will provide the detailed proof of Theorem~\ref{thm:main}.
\subsection{Analysis Related to Nesterov's Acceleration }

Because the update form of gradient tracking is almost the same to the one of SARAH \citep{nguyen2017sarah} which is a variance reduction method,  our main proof follows from proof framework of \cite{driggs2021accelerating} which is used to prove the convergence rate of accelerated variance reduction methods with proximal mapping. 

First, we give the relation between average variables in the following lemma.

\begin{lemma}
	Letting $\bbx_t$, $\bby_t$ and $\bbz_t$ (refer to Eq.~\eqref{eq:bbx}) be the average variables of $\xb_t$, $\yb_t$, $\zb_t$ respectively, then it holds that 
	\begin{align}
		\bbx_{t+1} =& \tau_t \bbz_t + (1 - \tau_t)\bby_t \label{eq:bbx_up}\\ 
		\bby_{t+1} =& \tau_t \bbz_{t+1} + (1 -\tau_t) \bby_t. \label{eq:yzy}
	\end{align}
\end{lemma}
\begin{proof}
	By multiplying $\frac{1}{m}\mathbf{1}^\top$ both sides of Eq.~\eqref{eq:xzy} and~\eqref{eq:yzy1} respectively and using the definition of average variables defined in Eq.~\eqref{eq:bbx}, we can obtain the result.
\end{proof}

Next, under Assumption~\ref{ass:f}, we give two properties similar to the ones of convexity and smoothness with the inexact gradient $\bbs_t$. 
\begin{lemma}[Lemma 1 of \cite{li2021accelerated}]
	Define 
	\begin{equation}\label{eq:ff}
		f(\bbx_t, \xb_t) \triangleq \frac{1}{m} \sum_{i=1}^{m} \left( f_i(\xbti) + \dotprod{\nabla f_i(\xbti), \bbx_t - \xbti} \right)
	\end{equation}
	Supposing Assumption~\ref{ass:f} holds, then we have for any $w \in \RR^d$,
	\begin{align}
		f(w) \ge& f(\bbx_t, \xb_t) + \dotprod{\bbs_t, w - \bbx_t},    \label{eq:cvx}\\
		f(w) \leq& f(\bbx_t, \xb_t) + \dotprod{\bbs_t, w - \bbx_t} + \frac{L}{2} \norm{w - \bbx_t}^2 + \frac{L}{2m} \norm{\Pi \xb_t}^2.  \label{eq:L}
	\end{align}
\end{lemma}

Next, applying the coupling framework, we construct a lower bound on the one-iteration progress of Algorithm~\ref{alg:alg_name}.

\begin{lemma}
	Letting Assumption~\ref{ass:f} hold, sequences $\{\xb_t\}$, $\{\yb_t\}$, $\{\zb_t\}$ and $\{\bs_t\}$ generated by Algorithm~\ref{alg:alg_name} satisfy the following property
	\begin{equation}\label{eq:main1}
		\begin{aligned}
			\gamma_t\left(f(\bbx_{t+1}, \xb_{t+1}) - f(x^*)\right)
			\leq& 
			\frac{\gamma_t(1-\tau_t)}{\tau_t} \Big( f(\bby_t) - f(\bbx_{t+1}, \xb_{t+1}) \Big) 
			- \frac{\gamma_t(1-\tau_t)}{\tau_t} D_f(\bby_t, \xb_{t+1}) \\
			&+\frac{\gamma_t}{\tau_t} \dotprod{\bbs_{t+1}, \bbx_{t+1} - \bby_{t+1} } + \gamma_t \dotprod{\bbs_{t+1}, \bbz_{t+1} - x^*},
		\end{aligned}
	\end{equation}
	where $f(\bbx_{t+1}, \xb_{t+1})$ and $D_f(\bby_t, \xb_{t+1})$ are defined in Eq.~\eqref{eq:ff} and Eq.~\eqref{eq:df}, respectively.
\end{lemma}
\begin{proof}
	By the convexity of $f_i(\cdot)$ and update rule of Algorithm~\ref{alg:alg_name}, we can obtain
	\begin{align*}
		&\gamma_t\left(f(\bbx_{t+1}, \xb_{t+1}) - f(x^*)\right)
		\stackrel{\eqref{eq:cvx}}{\leq}
		\gamma_t \dotprod{\bbs_{t+1}, \bbx_{t+1} - x^*}\\
		=&
		\gamma_t \dotprod{\bbs_{t+1}, \bbx_{t+1} - \bbz_t  } + \gamma_t\dotprod{\bbs_{t+1}, \bbz_t - x^* }\\
		\stackrel{\eqref{eq:bbx_up}}{=}&
		\frac{\gamma_t(1-\tau_t)}{\tau_t} \dotprod{\bbs_{t+1}, \bby_t - \bbx_{t+1}} + \gamma_t \dotprod{\bbs_{t+1}, \bbz_t - x^* }\\
		\stackrel{\eqref{eq:ff}\eqref{eq:df}}{=}&
		\frac{\gamma_t(1-\tau_t)}{\tau_t} \Big( f(\bby_t) - f(\bbx_{t+1}, \xb_{t+1}) - D_f(\bby_t, \xb_{t+1}) \Big) 
		+ \gamma_t \dotprod{\bbs_{t+1}, \bbz_t - x^* }\\
		=&
		\frac{\gamma_t(1-\tau_t)}{\tau_t} \Big( f(\bby_t) - f(\bbx_{t+1}, \xb_{t+1}) \Big) 
		- \frac{\gamma_t(1-\tau_t)}{\tau_t} D_f(\bby_t, \xb_{t+1}) \\
		&+ \gamma_t \dotprod{\bbs_{t+1}, \bbz_t - \bbz_{t+1} } + \gamma_t \dotprod{\bbs_{t+1}, \bbz_{t+1} - x^*}\\
		\stackrel{\eqref{eq:bbx_up}\eqref{eq:yzy}}{=}&
		\frac{\gamma_t(1-\tau_t)}{\tau_t} \Big( f(\bby_t) - f(\bbx_{t+1}, \xb_{t+1}) \Big) 
		- \frac{\gamma_t(1-\tau_t)}{\tau_t} D_f(\bby_t, \xb_{t+1}) \\
		&+\frac{\gamma_t}{\tau_t} \dotprod{\bbs_{t+1}, \bbx_{t+1} - \bby_{t+1} } + \gamma_t \dotprod{\bbs_{t+1}, \bbz_{t+1} - x^*}.
	\end{align*}
\end{proof}

Next, we bound the value of $ \frac{\gamma_t}{\tau_t} \dotprod{\bbs_{t+1}, \bbx_{t+1} - \bby_{t+1} } $ and $\gamma_t \dotprod{\bbs_{t+1}, \bbz_{t+1} - x^*}$ in the next two lemmas.

\begin{lemma}
	Letting Assumption~\ref{ass:g} hold, that is, $g(\cdot)$ is $\mu$-strongly convex and $\hzb_{t+1}$ be defined in Eq.~\eqref{eq:hz}, then  it holds that
	\begin{equation}\label{eq:gg}
		\begin{aligned}
			\gamma_t\dotprod{\bbs_{t+1}, \bbz_{t+1} - x^*}
			\leq&
			\frac{1}{2m} \norm{\zb_t - \mathbf{1}x^*}^2 - \frac{1+\mu\gamma_t}{2m} \norm{\hzb_{t+1} - \mathbf{1}x^*}^2 - \frac{1}{2\tau_t^2} \norm{\bbx_{t+1} - \bby_{t+1}}^2 + \gamma_t g(x^*) \\
			&- \gamma_t \cdot \frac{1}{m}\sum_{i=1}^{m}g(\hzb_{t+1}^{(i)})  
			+ \frac{1}{2m}\left( \gamma_t^2\norm{\Pi\bs_{t+1}}^2 + \norm{\Pi\hzb_{t+1}}^2 \right).
		\end{aligned}
	\end{equation}
\end{lemma}
\begin{proof}
	By Lemma~\ref{lem:g} with $z = \hzb_{t+1}^{(i)}$, $x = \zbti$, $d = \bs_{t+1}^{(i)}$, and $y = x^*$, we have
	\begin{align*}
		\gamma_t \dotprod{\bs_{t+1}^{(i)}, \hzb_{t+1}^{(i)} - x^*} 
		\leq
		\frac{1}{2} \norm{\zbti - x^*}^2 - \frac{1+\mu\gamma_t}{2} \norm{\hzb_{t+1}^{(i)} - x^*} - \frac{1}{2} \norm{\hzb_{t+1}^{(i)} - \zbti}^2 - \gamma_t g(\hzb_{t+1}^{(i)}) + \gamma_t g(x^*). 
	\end{align*}
	We also have
	\begin{align*}
		&\gamma_t\dotprod{\bbs_{t+1}, \bbz_{t+1} - x^*} \\
		=&
		\gamma_t \cdot  \frac{1}{m}\sum_{i=1}^{m} \dotprod{\bs_{t+1}^{(i)}, \hzb_{t+1}^{(i)} - x^*}  
		+ 
		\gamma_t\left(\dotprod{\bbs_{t+1}, \bbz_{t+1} - x^*} -  \frac{1}{m}\sum_{i=1}^{m} \dotprod{\bs_{t+1}^{(i)}, \hzb_{t+1}^{(i)} - x^*}\right)\\
		=&
		\gamma_t \cdot  \frac{1}{m}\sum_{i=1}^{m} \dotprod{\bs_{t+1}^{(i)}, \hzb_{t+1}^{(i)} - x^*} + \gamma_t \cdot \frac{1}{m}\sum_{i=1}^{m} \dotprod{\bs_{t+1}^{(i)},  \bbz_{t+1} -\hzb_{t+1}^{(i)}}\\
		=& 
		\gamma_t \cdot  \frac{1}{m}\sum_{i=1}^{m} \dotprod{\bs_{t+1}^{(i)}, \hzb_{t+1}^{(i)} - x^*} + \gamma_t \cdot \frac{1}{m}\sum_{i=1}^{m} \dotprod{\bs_{t+1}^{(i)} - \bbs_{t+1}, \bbz_{t+1} -\hzb_{t+1}^{(i)}}\\
		\stackrel{\eqref{eq:cauchy}}{\leq}&
		\gamma_t \cdot  \frac{1}{m}\sum_{i=1}^{m} \dotprod{\bs_{t+1}^{(i)}, \hzb_{t+1}^{(i)} - x^*} +   \frac{1}{m}\sum_{i=1}^{m} \left(\frac{\gamma_t^2\norm{\bs_{t+1}^{(i)} - \bbs_{t+1}}^2 + \norm{\hzb_{t+1}^{(i)} - \bbz_{t+1}}^2}{2}\right) \\
		=&
		\gamma_t \cdot  \frac{1}{m}\sum_{i=1}^{m} \dotprod{\bs_{t+1}^{(i)}, \hzb_{t+1}^{(i)} - x^*} + \frac{1}{2m}\left( \gamma_t^2\norm{\Pi\bs_{t+1}}^2 + \norm{\Pi\hzb_{t+1}}^2 \right),
	\end{align*}
	where the third equality is because of 
	\begin{align*}
		\sum_{i=1}^{m}\dotprod{\bbs_{t+1}, \hzb_{t+1}^{(i)} - \bbz_{t+1}} 
		= \dotprod{\bbs_{t+1}, \sum_{i=1}^{m}\hzb_{t+1}^{(i)} - m\bbz_{t+1}} 
		= \dotprod{\bbs_{t+1}, \sum_{i=1}^{m}\zb_{t+1}^{(i)} - m\bbz_{t+1}} 
		= 0.
	\end{align*}
	Using above two equations, we can obtain that
	\begin{align*}
		\gamma_t\dotprod{\bbs_{t+1}, \bbz_{t+1} - x^*}
		\leq&
		\frac{1}{2m} \norm{\zb_t - \mathbf{1}x^*}^2 - \frac{1+\mu\gamma_t}{2m} \norm{\hzb_{t+1} - \mathbf{1}x^*}^2 - \frac{1}{2m} \norm{\hzb_{t+1} - \zb_t}^2 \\
		&- \gamma_t \cdot \frac{1}{m}\sum_{i=1}^{m}g(\hzb_{t+1}^{(i)}) + \gamma_t g(x^*) 
		+ \frac{1}{2m}\left( \gamma_t^2\norm{\Pi\bs_{t+1}}^2 + \norm{\Pi\hzb_{t+1}}^2 \right).
	\end{align*}
	Furthermore,
	\begin{align*}
		&-\norm{\hzb_{t+1} - \zb_t}^2 \\
		=& 
		-\norm{  \hzb_{t+1} - \mathbf{1}\bbz_{t+1} - (\zb_t - \mathbf{1} \bbz_t) + (\mathbf{1} \bbz_{t+1} - \mathbf{1} \bbz_t)}^2\\
		=&
		-\left(\norm{ \Pi \hzb_{t+1} }^2 + \norm{\Pi \zb_t}^2 + m\norm{\bbz_{t+1} - \bbz_t}^2\right)  -2 \sum_{i=1}^{m}\dotprod{\hzb_{t+1}^{(i)} - \bbz_{t+1} - (\zbti -  \bbz_t), \bbz_{t+1} - \bbz_t} \\
		&+2\sum_{i=1}^{m}\dotprod{\hzb_{t+1}^{(i)} - \bbz_{t+1}, \zbti - \bbz_t}\\
		=&
		-\left(\norm{ \Pi \hzb_{t+1} }^2 + \norm{\Pi \zb_t}^2 + m\norm{\bbz_{t+1} - \bbz_t}^2\right)  + 2\sum_{i=1}^{m}\dotprod{\hzb_{t+1}^{(i)} - \bbz_{t+1}, \zbti - \bbz_t}\\
		\stackrel{\eqref{eq:cauchy}}{\leq}&
		-\left(\norm{ \Pi \hzb_{t+1} }^2 + \norm{\Pi \zb_t}^2 + m\norm{\bbz_{t+1} - \bbz_t}^2\right)  + \sum_{i=1}^{m} \left(\norm{\hzb_{t+1}^{(i)} - \bbz_{t+1}}^2 + \norm{\zbti - \bbz_t }^2\right)\\
		=&
		-m\norm{\bbz_{t+1} - \bbz_t}^2  
		\stackrel{\eqref{eq:bbx_up}\eqref{eq:yzy}}{=}
		-\frac{m}{\tau_t^2} \norm{\bbx_{t+1} - \bby_{t+1}}^2,
	\end{align*}
	where the third equality is because of 
	\begin{align*}
		\sum_{i=1}^{m}\dotprod{\hzb_{t+1}^{(i)} - \bbz_{t+1} - (\zbti -  \bbz_t), \bbz_{t+1} - \bbz_t} 
		= 
		\dotprod{\sum_{i=1}^{m} \hzb_{t+1}^{(i)} - m \bbz_{t+1} +  \sum_{i=1}^{m} \zbti - m\bbz_t,  \bbz_{t+1} - \bbz_t} = 0.
	\end{align*}
	
	Therefore, we can obtain that
	\begin{align*}
		\gamma_t\dotprod{\bbs_{t+1}, \bbz_{t+1} - x^*}
		\leq&
		\frac{1}{2m} \norm{\zb_t - \mathbf{1}x^*}^2 - \frac{1+\mu\gamma_t}{2m} \norm{\hzb_{t+1} - \mathbf{1}x^*}^2 - \frac{1}{2\tau_t^2} \norm{\bbx_{t+1} - \bby_{t+1}}^2 + \gamma_t g(x^*) \\
		&- \gamma_t \cdot \frac{1}{m}\sum_{i=1}^{m}g(\hzb_{t+1}^{(i)})  
		+ \frac{1}{2m}\left( \gamma_t^2\norm{\Pi\bs_{t+1}}^2 + \norm{\Pi\hzb_{t+1}}^2 \right).
	\end{align*}
\end{proof}

\begin{lemma}
	Letting Assumption~\ref{ass:f} and Assumption~\ref{ass:g} hold, and $\hzb_{t+1}$ be defined in Eq.~\eqref{eq:hz}, then it holds that
	\begin{equation}\label{eq:sxy}
		\begin{aligned}
			\frac{\gamma_t}{\tau_t} \dotprod{ \bbs_{t+1}, \bbx_{t+1} - \bby_{t+1} } 
			\leq& 
			\frac{\gamma_t}{\tau_t} \Big(f(\bbx_{t+1}, \xb_{t+1}) - F(\bby_{t+1})\Big) +   \frac{\gamma_t}{m} \sum_{i=1}^{m} g(\hzb_{t+1}^{(i)}) + \frac{\gamma_t (1-\tau_t)}{\tau_t} g(\bby_t) \\
			&+ \frac{L\gamma_t}{2\tau_t}\norm{\bby_{t+1} - \bbx_{t+1}}^2 + \frac{L\gamma_t}{2m\tau_t} \norm{\Pi \xb_{t+1}}^2.
		\end{aligned}
	\end{equation}
\end{lemma}
\begin{proof}
	By the $L$-smoothness of $f_i(\cdot)$ and the convexity of $g(\cdot)$,  we have
	\begin{align*}
		&\frac{\gamma_t}{\tau_t} \dotprod{ \bbs_{t+1}, \bbx_{t+1} - \bby_{t+1} } \\
		\stackrel{\eqref{eq:L}}{\leq}&
		\frac{\gamma_t}{\tau_t} \Big(f(\bbx_{t+1}, \xb_{t+1}) - f(\bby_{t+1})\Big) + \frac{L\gamma_t}{2\tau_t}\norm{\bby_{t+1} - \bbx_{t+1}}^2 + \frac{L\gamma_t}{2m\tau_t} \norm{\Pi \xb_{t+1}}^2\\
		=& 
		\frac{\gamma_t}{\tau_t} \Big(f(\bbx_{t+1}, \xb_{t+1}) - F(\bby_{t+1})\Big) + \frac{\gamma_t}{\tau_t} g(\bby_{t+1}) + \frac{L\gamma_t}{2\tau_t}\norm{\bby_{t+1} - \bbx_{t+1}}^2 + \frac{L\gamma_t}{2m\tau_t} \norm{\Pi \xb_{t+1}}^2\\
		\leq&
		\frac{\gamma_t}{\tau_t} \Big(f(\bbx_{t+1}, \xb_{t+1}) - F(\bby_{t+1})\Big) + \gamma_t g(\bbz_{t+1}) + \frac{\gamma_t (1-\tau_t)}{\tau_t} g(\bby_t) + \frac{L\gamma_t}{2\tau_t}\norm{\bby_{t+1} - \bbx_{t+1}}^2 + \frac{L\gamma_t}{2m\tau_t} \norm{\Pi \xb_{t+1}}^2\\
		\leq&
		\frac{\gamma_t}{\tau_t} \Big(f(\bbx_{t+1}, \xb_{t+1}) - F(\bby_{t+1})\Big) + \frac{\gamma_t}{m} \sum_{i=1}^{m} g(\hzb_{t+1}^{(i)}) + \frac{\gamma_t (1-\tau_t)}{\tau_t} g(\bby_t) + \frac{L\gamma_t}{2\tau_t}\norm{\bby_{t+1} - \bbx_{t+1}}^2 + \frac{L\gamma_t}{2m\tau_t} \norm{\Pi \xb_{t+1}}^2,
	\end{align*}
	where the second inequality is because of Eq.~\eqref{eq:yzy} and the convexity of $g(\cdot)$.
	The last inequality is because of $\bbz_{t+1} = \frac{1}{m}\sum_{i=1}^{m}\zb_{t+1}^{(i)} = \frac{1}{m}\sum_{i=1}^{m}\hzb_{t+1}^{(i)}$ and the convexity of $g(\cdot)$.
\end{proof}

\begin{lemma}\label{lem:main}
	Letting Assumption~\ref{ass:f} and Assumption~\ref{ass:g} hold, and $\hzb_{t+1}$ be defined in Eq.~\eqref{eq:hz}, then it holds that
	\begin{equation}\label{eq:main2}
		\begin{aligned}
			0 \leq& \frac{\gamma_t(1-\tau_t)}{\tau_t} F(\bby_t) - \frac{\gamma_t}{\tau_t}F(\bby_{t+1}) + \gamma_t F(x^*) + \left(\frac{L\gamma_t}{2\tau_t} - \frac{1}{2\tau_t^2}\right) \norm{\bbx_{t+1} - \bby_{t+1}}^2 \\
			&+\frac{1}{2m} \norm{\zb_t - \mathbf{1}x^*}^2 - \frac{1+\mu\gamma_t}{2m} \norm{\hzb_{t+1} - \mathbf{1}x^*}^2 - \frac{\gamma_t(1-\tau_t)}{\tau_t} D_f(\bby_t, \xb_{t+1})\\
			&+ \frac{L\gamma_t}{2m\tau_t} \norm{\Pi \xb_{t+1}}^2 + \frac{1}{2m}\left( \gamma_t^2\norm{\Pi\bs_{t+1}}^2 + \norm{\Pi\hzb_{t+1}}^2 \right).
		\end{aligned}
	\end{equation}
\end{lemma}
\begin{proof}
	Combining Eq.~\eqref{eq:main1}, Eq.~\eqref{eq:gg} and Eq.~\eqref{eq:sxy}, we can obtain that
	\begin{align*}
		0 \leq& - \gamma_t \Big(f(\bbx_{t+1}, \xb_{t+1}) - f(x^*)\Big) +  \frac{\gamma_t(1-\tau_t)}{\tau_t} \Big( f(\bby_t) - f(\bbx_{t+1}, \xb_{t+1}) \Big) - \frac{\gamma_t(1-\tau_t)}{\tau_t} D_f(\bby_t, \xb_{t+1})\\
		&+\frac{1}{2m} \norm{\zb_t - \mathbf{1}x^*}^2 - \frac{1+\mu\gamma_t}{2} \norm{\hzb_{t+1} - \mathbf{1}x^*}^2 - \frac{1}{2\tau_t^2} \norm{\bbx_{t+1} - \bby_{t+1}}^2 + \gamma_t g(x^*) \\
		&- \gamma_t \cdot \frac{1}{m}\sum_{i=1}^{m}g(\hzb_{t+1}^{(i)})  
		+ \frac{1}{2m}\left( \gamma_t^2\norm{\Pi\bs_{t+1}}^2 + \norm{\Pi\hzb_{t+1}}^2 \right)\\
		&+\frac{\gamma_t}{\tau_t} \Big(f(\bbx_{t+1}, \xb_{t+1}) - F(\bby_{t+1})\Big) +   \frac{\gamma_t}{m} \sum_{i=1}^{m} g(\hzb_{t+1}^{(i)}) + \frac{\gamma_t (1-\tau_t)}{\tau_t} g(\bby_t) \\
		&+ \frac{L\gamma_t}{2\tau_t}\norm{\bby_{t+1} - \bbx_{t+1}}^2 + \frac{L\gamma_t}{2m\tau_t} \norm{\Pi \xb_{t+1}}^2\\
		=&
		\frac{\gamma_t(1-\tau_t)}{\tau_t} F(\bby_t) - \frac{\gamma_t}{\tau_t}F(\bby_{t+1}) + \gamma_t F(x^*) + \left(\frac{L\gamma_t}{2\tau_t} - \frac{1}{2\tau_t^2}\right) \norm{\bbx_{t+1} - \bby_{t+1}}^2 \\
		&+\frac{1}{2m} \norm{\zb_t - \mathbf{1}x^*}^2 - \frac{1+\mu\gamma_t}{2} \norm{\hzb_{t+1} - \mathbf{1}x^*}^2 - \frac{\gamma_t(1-\tau_t)}{\tau_t} D_f(\bby_t, \xb_{t+1})\\
		&+ \frac{L\gamma_t}{2m\tau_t} \norm{\Pi \xb_{t+1}}^2 + \frac{1}{2m}\left( \gamma_t^2\norm{\Pi\bs_{t+1}}^2 + \norm{\Pi\hzb_{t+1}}^2 \right).
	\end{align*}
\end{proof}
\begin{remark}
	Lemma~\ref{lem:main} shares almost the same result as the one of Lemma~4 of \cite{driggs2021accelerating}. However, in our lemma, we have extra consensus terms $\norm{ \Pi\xb_{t+1} }^2$, $\norm{\Pi\hzb_{t+1}}^2$ and $ \norm{\bs_{t+1}}^2$ which are because of the decentralized optimization setting.
	Furthermore, an extended Bregman divergence $D_f(\bby_t, \xb_{t+1})$ is used in our lemma while the standard  Bregman divergence is used in Lemma~4 of \cite{driggs2021accelerating}.	
\end{remark}
\begin{remark}
	Our algorithm is almost the same to the one of \cite{li2021accelerated} except the difference in the proximal mapping in the Eq.~\eqref{eq:prox_step}. 
	However, the proof framework of \cite{li2021accelerated} can not be extended to the composite optimization. The proof in \citep{li2021accelerated} relies the following identity $$\tau \dotprod{\bbs_{t+1}, \bbz_{t+1} - x^*} = \frac{\tau^2}{2\gamma}\left( \norm{ \bbz_t - x^*  }^2 - \norm{\bbz_{t+1} - x^*}^2 - \norm{\bbz_{t+1} - \bbz_t}^2 \right) + \dots$$
	which implies that
	\begin{align*}
		\gamma \dotprod{\bbs_{t+1}, \bbz_{t+1} - x^*} = \frac{\tau}{2} \left( \norm{ \bbz_t - x^*  }^2 - \norm{\bbz_{t+1} - x^*}^2 - \norm{\bbz_{t+1} - \bbz_t}^2 \right) + \dots
	\end{align*}
	Comparing above equation with Eq.~\eqref{eq:gg}, we can observe that above equation scaling Eq.~\eqref{eq:gg} with a factor $\tau$. 
	Due to this difference, the  proof framework of \cite{li2021accelerated} can not be extended to prove the convergence rate of our algorithm.
\end{remark}

\subsection{Consensus Error Bound}

This section, we will bound the consensus error terms in Lemma~\ref{lem:main}. Parts of proofs in this section follow from the ones of \citet{li2021accelerated}.
First, we prove a property related to the proximal mapping in the decentralized setting. 
\begin{lemma}
	For any $\xb\in\RR^{m\times d}$, the proximal mapping $\prox_{\gamma g}$ has the following property,
	\begin{equation}\label{eq:prox_con}
		\norm{\prox_{\gamma g} \left(\frac{1}{m}\mathbf{1}\mathbf{1}^\top\xb\right) - \frac{1}{m}\mathbf{1}\mathbf{1}^\top \prox_{\gamma g}(\xb)}
		\leq 
		\norm{\Pi\xb}.
	\end{equation}
\end{lemma}
\begin{proof}
	Using the non-expansiveness of the proximal mapping, we have
	\begin{align*}
		&\norm{\prox_{\gamma g} \left(\frac{1}{m}\mathbf{1}\mathbf{1}^\top\xb\right) - \frac{1}{m}\mathbf{1}\mathbf{1}^\top \prox_{\gamma g}(\xb)}^2
		=m \norm{\prox_{\gamma g}\left(\frac{1}{m} \mathbf{1}^\top \xb\right) -\frac{1}{m} \sum_{i=1}^{m} \prox_{\gamma g}(\xb^{(i)}) }^2\\
		=&
		m\norm{\frac{1}{m}\sum_{i=1}^{m} \left(\prox_{\gamma g} \left(\frac{1}{m}\mathbf{1}^\top \xb \right) - \prox_{\gamma g}(\xb^{(i)})\right)}^2
		\leq
		\sum_{i=1}^{m} \norm{ \prox_{\gamma g} \left(\frac{1}{m}\mathbf{1}^\top \xb \right) - \prox_{\gamma g}(\xb^{(i)}) }^2 \\
		\leq&
		\sum_{i=1}^{m}
		\norm{ \frac{1}{m}\mathbf{1}^\top \xb - \xb^{(i)}}^2 = \norm{\xb - \frac{1}{m}\mathbf{1}\mathbf{1}^\top \xb}^2 = \norm{\Pi\xb}^2.
	\end{align*}
\end{proof}

Next, we will give the upper bounds of consensus error terms and their convergence properties. 
\begin{lemma}\label{lem:err}
	Letting $K = \frac{15}{\sqrt{1 - \lambda_2(W)}}$, and the step size $\gamma_t$ be properly chosen such that $L^2\gamma_t^2 \leq \frac{1}{100\cdot \tau_t^2}$ in Algorithm~\ref{alg:alg_name}, sequences $\{\xb_t\}$, $\{\yb_t\}$, $\{\hzb_t\}$ (defined in Eq.~\eqref{eq:hz})  and $\{\bs_t\}$ satisfy that
	\begin{equation}\label{eq:err}
		\begin{aligned}
			\norm{\Pi \xb_{t+1}}^2 \leq& 2\norm{\Pi\yb_t}^2 + \frac{\tau_t^2}{32}\norm{\Pi\hzb_t}^2,\qquad
			\norm{\Pi\hzb_{t+1}}^2 \leq  8 \left( \frac{1}{64} \norm{\Pi\hzb_t}^2  + \gamma_t^2 \norm{\Pi\bs_{t+1}}^2 \right),\\
			\gamma_t^2\norm{\Pi\bs_{t+1}}^2 \leq& \frac{1}{32} \cdot \gamma_t^2\norm{\Pi\bs_t  }^2 + 8\gamma_t^2\left(mL D_f (\bby_t, \xb_{t+1}) + mL^2 \norm{\bbx_t - \bby_t}^2 + L^2 \norm{\Pi\xb_t}^2\right),
		\end{aligned}
	\end{equation}
	and
	\begin{equation}\label{eq:xyzs}
		\begin{aligned}
			&\norm{\Pi\xb_{t+1}}^2 + c_1 \norm{\Pi \yb_{t+1}}^2 + c_2\norm{\Pi\hzb_{t+1}}^2 + c_3  \gamma_t^2 \norm{\Pi\bs_{t+1}}^2\\
			\leq& 
			\frac{1}{4}\left( \norm{\Pi\xb_t}^2
			+
			c_1\norm{\Pi\yb_t}^2 +  c_2\norm{\Pi\hzb_t}^2 +  c_3\gamma_t^2\norm{\Pi\hzb_t}^2 \right)\\
			&+25\tau_t^2\gamma_t^2\left(mL D_f (\bby_t, \xb_{t+1}) + mL^2 \norm{\bbx_t - \bby_t}^2\right),
		\end{aligned}
	\end{equation}
	with $c_1 = 64$, $c_2 = \frac{9\tau_t^2}{32}$, and $c_3 = \frac{5\tau_t^2}{14}$.
\end{lemma}
\begin{proof}
	For the notation convenience, we denote that
	\begin{align*}
		\TB(\xb) \triangleq \fm(\xb, K).
	\end{align*}
	By the chosen communication step number $K$ in ``FastMix'' and Proposition~\ref{lem:mix_eq}, we can obtain that
	\begin{align*}
		\norm{\Pi\cdot \TB(\xb)} \leq \rho \norm{ \Pi \xb }, \mbox{ with } \rho^2 \leq \frac{1}{64}.
	\end{align*}
	First, by  the definition of $\hzb_{t+1}$ in Eq.~\eqref{eq:hz} and the update rule of $\zb_{t+1}$ in Eq~\eqref{eq:prox_step}, we can obtain
	\begin{align*}
		\norm{\Pi \hzb_{t+1}} =  \norm{ \Pi \cdot \prox_{\gamma_t g}(\zb_t - \gamma_t \bs_{t+1}) }.
	\end{align*}
	Furthermore,
	\begin{align*}
		&\norm{ \Pi \cdot \prox_{\gamma_t g}(\zb_t - \gamma_t \bs_{t+1}) } 
		=
		\norm{  \prox_{\gamma_t g}(\zb_t - \gamma_t \bs_{t+1}) - \frac{1}{m}\mathbf{1}\mathbf{1}^\top \prox_{\gamma_t g}(\zb_t - \gamma_t \bs_{t+1}) }  \\
		\leq&
		\norm{ \prox_{\gamma_t g}(\zb_t - \gamma_t \bs_{t+1}) - \prox_{\gamma_t g}(\mathbf{1}(\bbz_t - \gamma_t \bbs_{t+1}))  } \\
		&+
		\norm{\prox_{\gamma_t g}(\mathbf{1}(\bbz_t - \gamma_t \bbs_{t+1})) - \frac{1}{m}\mathbf{1}\mathbf{1}^\top \prox_{\gamma_t g}(\zb_t - \gamma_t \bs_{t+1})}\\
		\leq&
		\norm{\zb_t - \mathbf{1}\bbz_t} + \gamma_t\norm{ \bs_{t+1} - \mathbf{1}\bbs_{t+1} } +
		\norm{\prox_{\gamma_t g}(\mathbf{1}(\bbz_t - \gamma_t \bbs_{t+1})) - \frac{1}{m}\mathbf{1}\mathbf{1}^\top \prox_{\gamma_t g}(\zb_t - \gamma_t \bs_{t+1})}\\
		\stackrel{\eqref{eq:prox_con}}{\leq}&
		\norm{\Pi \zb_t} + \gamma_t\norm{\Pi\bs_{t+1}} + \norm{\Pi (\zb_t - \gamma_t \bs_{t+1})} 
		\leq
		2 \norm{ \Pi\zb_t } + 2\gamma_t \norm{\Pi \bs_{t+1}}\\
		\stackrel{\eqref{eq:prox_step}}{\leq}&
		2\rho\norm{\Pi\hzb_t} + 2\gamma_t \norm{\Pi \bs_{t+1}},
	\end{align*}
	where the second inequality is because of the non-expansiveness of the proximal mapping and the last inequality is because of $\zb_t = \TB(\hzb_t)$.
	Thus, it holds that
	\begin{align}
		\norm{\Pi\hzb_{t+1}}^2 \leq  \big( 2\rho \norm{ \Pi\hzb_t } + 2\gamma_t \norm{\Pi \bs_{t+1}} \big)^2 \leq  8 \left(\rho^2 \norm{\Pi\hzb_t}^2  + \gamma_t^2 \norm{\Pi\bs_{t+1}}^2 \right). \label{eq:zz}
	\end{align}
	Furthermore, by the update rule of $\xb_{t+1}$ in Eq.~\eqref{eq:xzy}, we have
	\begin{align*}
		\norm{\Pi \xb_{t+1}}^2
		\leq& 
		\left(\tau_t \norm{\Pi \zb_t} + (1 - \tau_t) \norm{\Pi \yb_t} \right)^2
		\leq 
		2\tau_t^2\norm{\Pi\zb_t}^2 + 2\norm{\Pi\yb_t}^2
		\stackrel{\eqref{eq:prox_step}}{\leq}
		2\rho^2\tau_t^2\norm{\Pi\hzb_t}^2 + 2\norm{\Pi\yb_t}^2.
	\end{align*}
	Similarly, by the update rule of $\yb_{t+1}$ in Eq.~\eqref{eq:yzy1}, we have
	\begin{align*}
	\norm{\Pi \yb_{t+1}}^2 
		\leq& 
		\left(\tau_t  \rho\norm{\Pi\zb_{t+1}} + (1-\tau_t)\rho \norm{\Pi \yb_t}\right)^2 
		\leq
		2\rho^2\norm{\Pi\yb_t}^2 + 2\tau_t^2\rho^4 \norm{\Pi\hzb_{t+1}}^2
	\end{align*}
	
	By the update rule of $\bs_{t+1}$ in Eq.~\eqref{eq:ss}, we have
	\begin{align*}
		&\norm{ \Pi \bs_{t+1} }^2 
		= 
		\norm{\Pi\cdot\TB\big( \bs_t + \nabla \widetilde{f}(\xb_{t+1}) - \nabla \widetilde{f}(\xb_t)\big)}^2\\
		\leq&
		\rho^2 \norm{\Pi\left( \bs_t + \nabla \widetilde{f}(\xb_{t+1}) - \nabla \widetilde{f}(\xb_t) \right)}^2
		\leq
		2\rho^2 \norm{ \Pi \bs_t }^2 + 2\norm{\nabla \widetilde{f}(\xb_{t+1}) - \nabla \widetilde{f}(\xb_t)  }^2.
	\end{align*}
	Furthermore, 
	\begin{align*}
		&\norm{\nabla \widetilde{f}(\xb_{t+1}) - \nabla \widetilde{f}(\xb_t)}^2 
		=
		\sum_{i=1}^{m} \norm{\nabla f_i(\xb_{t+1}^{(i)}) - \nabla f_i(\xbti)}^2\\
		\stackrel{\eqref{eq:xy}}{\leq}&
		2\sum_{i=1}^{m} \norm{ \nabla f_i(\xb_{t+1}^{(i)}) - \nabla f_i(\bby_t) }^2 
		+ 
		4\sum_{i=1}^{m}   \left( \norm{ \nabla f_i (\bby_t) - \nabla f_i(\bbx_t) }^2 + \norm{ \nabla f_i(\bbx_t) - \nabla f_i(\xbti)}^2 \right)\\
		\stackrel{\eqref{eq:L_breg}\eqref{eq:df}}{\leq}&
		4mL D_f (\bby_t, \xb_{t+1})  
		+
		4\sum_{i=1}^{m}   \left( \norm{ \nabla f_i (\bby_t) - \nabla f_i(\bbx_t) }^2 + \norm{ \nabla f_i(\bbx_t) - \nabla f_i(\xbti)}^2 \right)\\
		\stackrel{\eqref{eq:L_n}}{\leq}&
		4mL D_f (\bby_t, \xb_{t+1})  
		+ 4L^2 \sum_{i=1}^{m} \left(\norm{\bbx_t - \bby_t}^2 + \norm{\xbti - \bbx_t}^2\right) \\
		=&
		4mL D_f (\bby_t, \xb_{t+1}) + 4mL^2 \norm{\bbx_t - \bby_t}^2 + 4L^2 \norm{\Pi\xb_t}^2.
	\end{align*}
	Thus, it holds that
	\begin{align}
		\gamma_t^2\norm{\Pi\bs_{t+1}}^2 \leq& 2\rho^2 \cdot \gamma_t^2\norm{\Pi\bs_t  }^2 + 8\gamma_t^2\left(mL D_f (\bby_t, \xb_{t+1}) + mL^2 \norm{\bbx_t - \bby_t}^2 + L^2 \norm{\Pi\xb_t}^2\right). \label{eq:sss}
	\end{align}
	Combining above results, we can obtain that
	\begin{align*}
	&\norm{\Pi\xb_{t+1}}^2 + c_1\norm{\Pi\yb_{t+1}}^2 + c_2\norm{\Pi\hzb_{t+1}}^2 + c_3\gamma_t^2\norm{\Pi\bs_{t+1}}^2\\
	\leq&
	2\rho^2\tau_t^2\norm{\Pi\hzb_t}^2 + 2\norm{\Pi\yb_t}^2 + 2c_1\rho^2\norm{\Pi\yb_t}^2 + 2c_1\tau_t^2\rho^4\norm{\Pi\hzb_{t+1}}^2 +c_2\norm{\Pi\hzb_{t+1}}^2 + c_3\gamma_t^2\norm{\Pi\bs_{t+1}}^2\\
	\stackrel{\eqref{eq:zz}}{\leq}&
	2(1+c_1\rho^2)\norm{\Pi\yb_t}^2 + 2\rho^2\tau_t^2\norm{\Pi\hzb_t}^2 + 8(2c_1\tau_t^2\rho^4 +c_2)\rho^2 \norm{\Pi\hzb_t}^2 + \big(8(2c_1\tau_t^2\rho^4 +c_2)+c_3\big) \gamma_t^2\norm{\Pi\bs_{t+1}}^2\\
	\stackrel{\eqref{eq:sss}}{\leq}&
	8\big(8(2c_1\tau_t^2\rho^4 +c_2)+c_3\big) L^2\gamma_t^2\norm{\Pi\xb_t}^2 + \frac{2(1+c_1\rho^2)}{c_1}\cdot c_1\norm{\Pi\yb_t}^2\\
	&+\frac{(8(2c_1\tau_t^2\rho^4 +c_2)+2\tau_t^2)\rho^2}{c_2}\cdot c_2\norm{\Pi\hzb_t}^2+\frac{2\big(8(2c_1\tau_t^2\rho^4 +c_2)+c_3\big)\rho^2}{c_3} \cdot c_3\gamma_t^2\norm{\Pi\bs_t}^2\\
	&+8\big(8(2c_1\tau_t^2\rho^4 +c_2)+c_3\big) \gamma_t^2\left(mL D_f (\bby_t, \xb_{t+1}) + mL^2 \norm{\bbx_t - \bby_t}^2\right)
	\end{align*}
	
	By the value of $c_1$, $c_2$ and $c_3$, we can obtain that
	\begin{align*}
	&\norm{\Pi\xb_{t+1}}^2 + c_1\norm{\Pi\yb_{t+1}}^2 + c_2\norm{\Pi\hzb_{t+1}}^2 + c_3\gamma_t^2\norm{\Pi\bs_{t+1}}^2\\
	\leq&
	25\tau_t^2 L^2\gamma_t^2 \norm{\Pi\xb_t}^2
	+
	\frac{1}{4} c_1\norm{\Pi\yb_t}^2 + \frac{1}{4} c_2\norm{\Pi\hzb_t}^2 + \frac{1}{4} c_3\gamma_t^2\norm{\Pi\hzb_t}^2 \\
	&+25\tau_t^2\gamma_t^2\left(mL D_f (\bby_t, \xb_{t+1}) + mL^2 \norm{\bbx_t - \bby_t}^2\right)\\
	\leq&
	\frac{1}{4}\left( \norm{\Pi\xb_t}^2
	+
	 c_1\norm{\Pi\yb_t}^2 +  c_2\norm{\Pi\hzb_t}^2 +  c_3\gamma_t^2\norm{\Pi\hzb_t}^2 \right)\\
	 &+25\tau_t^2\gamma_t^2\left(mL D_f (\bby_t, \xb_{t+1}) + mL^2 \norm{\bbx_t - \bby_t}^2\right).
	\end{align*}
\end{proof}

\subsection{Proof of Theorem~\ref{thm:main}}

Next, we will bound the consensus error terms in Eq.~\eqref{eq:main2} and sum them from $t = 0$ to $T-1$ scaled with a non-negative sequence $\{\beta_t\}$.
\begin{lemma}\label{lem:aa}
 Suppose Assumption~\ref{ass:f} and Assumption~\ref{ass:g} hold with $\mu > 0$.	
 Let $K = \frac{15}{\sqrt{1 - \lambda_2(W)}}$, and a constant step size $\gamma$ and $\tau$ are set as Theorem~\ref{thm:main}. 
 Denoting that $\beta_t = (1+\mu\gamma)^t$, then it holds that
	\begin{equation*}
		\begin{aligned}
			&\sum_{t=1}^{T-1} \beta_t \left(\frac{L\gamma}{2m\tau} \norm{\Pi \xb_{t+1}}^2 + \frac{1}{2m}\left( \gamma^2\norm{\Pi\bs_{t+1}}^2 + \norm{\Pi\hzb_{t+1}}^2  \right)\right)\\
			\leq&
			\frac{c_4}{2m}\sum_{t=1}^{T-1}\beta_t \left(\frac{1}{4}\right)^t\left(\norm{\Pi\xb_1}^2 + c_1\norm{\Pi\yb_1}^2 + c_2\norm{\Pi\hzb_1}^2 + c_3\gamma^2\norm{\Pi\bs_1}^2\right)\\ 
			&+ \frac{25c_4}{m}\sum_{t=1}^{T-1}\beta_t\tau^2\gamma^2\left(mL D_f (\bby_t, \xb_{t+1}) + mL^2 \norm{\bbx_t - \bby_t}^2\right). 
		\end{aligned}
	\end{equation*}
	where $c_4 \triangleq \max\left\{\frac{L\gamma}{\tau}, \frac{1}{c_2}, \frac{1}{c_3}\right\}$ and $c_1$, $c_2$, $c_3$ are defined in Lemma~\ref{lem:err}.
\end{lemma}
\begin{proof}
	First, by the setting of the step size $\gamma = \frac{1}{20\sqrt{L\mu}}$ and the definition $\tau =\mu\gamma $, we have
	\begin{align*}
		L^2\gamma^2 = \frac{L}{400\mu} \leq  \frac{1}{100\cdot \mu^2 \cdot \frac{1}{400\mu L}} = \frac{1}{100\cdot \tau^2}.
	\end{align*} 
	Thus, the condition $L^2\gamma^2 \leq\frac{1}{100\cdot \tau^2}$ in Lemma~\ref{lem:err} holds.
	Furthermore, by the setting of $K$ in Theorem~\ref{thm:main}, the results in Lemma~\ref{lem:err} hold.
	
	Let us denote that   $\Delta_t = 25\tau^2\gamma^2\left(mL D_f (\bby_t, \xb_{t+1}) + mL^2 \norm{\bbx_t - \bby_t}^2\right)$. 
	Then, by Lemma~\ref{lem:err}, 
	we can obtain that
	\begin{align*}
		&\sum_{t=1}^{T-1} \beta_t \left(\frac{L\gamma}{2m\tau} \norm{\Pi \xb_{t+1}}^2 + \frac{1}{2m}\left( \gamma^2\norm{\Pi\bs_{t+1}}^2 + \norm{\Pi\hzb_{t+1}}^2  \right)\right)\\
		\leq&
		\frac{1}{2m}\sum_{t=1}^{T-1} \beta_t \max\left\{\frac{L\gamma}{\tau}, \frac{1}{c_2}, \frac{1}{c_3}\right\} \left(\norm{\Pi\xb_{t+1}}^2 + c_1\norm{\Pi\yb_{t+1}}^2 + c_2\norm{\Pi\hzb_{t+1}}^2 + c_3\gamma^2\norm{\Pi\bs_{t+1}}^2\right)\\
		\leq&
		\frac{c_4}{2m}\sum_{t=1}^{T-1}\beta_t \left(\norm{\Pi\xb_{t+1}}^2 + c_1\norm{\Pi\yb_{t+1}}^2 + c_2\norm{\Pi\hzb_{t+1}}^2 + c_3\gamma^2\norm{\Pi\bs_{t+1}}^2\right)\\
		\leq&
		\frac{c_4}{2m}\sum_{t=1}^{T-1}\beta_t \left(\frac{1}{4}\right)^t\left(\norm{\Pi\xb_1}^2 + c_1\norm{\Pi\yb_1}^2 + c_2\norm{\Pi\hzb_1}^2 + c_3\gamma^2\norm{\Pi\bs_1}^2\right) 
		+ \frac{c_4}{2m}\sum_{t=1}^{T-1}\beta_t\sum_{k=1}^{t} \left(\frac{1}{4}\right)^{t-k}\Delta_k\\
		\leq&
		\frac{c_4}{2m}\sum_{t=1}^{T-1}\beta_t \left(\frac{1}{4}\right)^t\left(\norm{\Pi\xb_1}^2 + c_1\norm{\Pi\yb_1}^2 + c_2\norm{\Pi\hzb_1}^2 + c_3\gamma^2\norm{\Pi\bs_1}^2\right) 
		+ \frac{c_4}{2m}\sum_{t=1}^{T-1}\sum_{k=1}^{t}\beta_k \left(\frac{1}{2}\right)^{t-k}\Delta_k\\
		\leq&
		\frac{c_4}{2m}\sum_{t=1}^{T-1}\beta_t \left(\frac{1}{4}\right)^t\left(\norm{\Pi\xb_1}^2 + c_1\norm{\Pi\yb_1}^2 + c_2\norm{\Pi\hzb_1}^2 + c_3\gamma^2\norm{\Pi\bs_1}^2\right) 
		+ \frac{c_4}{m}\sum_{t=1}^{T-1}\beta_t\Delta_t,
	\end{align*}
	where the forth inequality is because of the fact
	\begin{equation*}
		\beta_t\cdot \frac{1}{4} = (1 + \mu\gamma)^t \cdot \frac{1}{4} = (1 + \mu\gamma)^{t-1} c\dot (1 + \mu\gamma)  \frac{1}{4} \leq \beta_{t-1}\cdot \frac{1}{2},
	\end{equation*}   
	and  the last inequality is because of Lemma~\ref{lem:recs} with $\rho = \frac{1}{2}$.
	
\end{proof}

Combining results in previous sections, we can prove Theorem~\ref{thm:main}.
\begin{proof}[Proof of Theorem~\ref{thm:main}]
	By the definition of  $c_1$, $c_2$, $c_3$, and $c_4$  which are defined in Lemma~\ref{lem:err} and Lemma~\ref{lem:aa} respectively, we can obtain that
	\begin{align*}
		c_4 =  \max\left\{\frac{L\gamma}{\tau}, \frac{1}{c_2}, \frac{1}{c_3}\right\} = \max\left\{ \frac{1}{400\tau^2}, \frac{32}{9\tau^2}, \frac{14}{5\tau^2} \right\}  \leq \frac{4}{\tau^2}.
	\end{align*}
	Then, the result of Lemma~\ref{lem:aa} can be represented as with $\beta_t = (1 + \mu\gamma)^t$, 
	\begin{equation}\label{eq:rec}
	\begin{aligned}
		&\sum_{t=0}^{T-1} \beta_t \left(\frac{L\gamma}{2m\tau} \norm{\Pi \xb_{t+1}}^2 + \frac{1}{2m}\left( \gamma^2\norm{\Pi\bs_{t+1}}^2 + \norm{\Pi\zb_{t+1}}^2 \right)\right)\\
	\leq&
	\frac{2}{m\tau^2}  \left( \norm{\Pi\xb_1}^2 + c_1 \norm{\Pi\yb_1}^2 + c_2 \norm{\Pi \zb_1}^2 + c_3\gamma^2 \norm{\Pi\bs_1}^2 \right) \cdot \sum_{t=0}^{T-1}\beta_t 4^{-t} \\
	&+ 100 \sum_{t=1}^{T-1} \beta_t \left( L\gamma^2 D_f(\bby_t, \xb_{t+1}) + L^2\gamma^2 \norm{\bbx_t - \bby_t}^2 \right)\\
	=&
	\frac{2}{m\tau^2}  \left( \norm{\Pi\xb_1}^2 + c_1 \norm{\Pi\yb_1}^2 + c_2 \norm{\Pi \zb_1}^2 + c_3\gamma^2 \norm{\Pi\bs_1}^2 \right) \cdot \sum_{t=0}^{T-1}\beta_t 4^{-t} \\
	&+ 100 \sum_{t=1}^{T-1} \beta_t  L\gamma^2 D_f(\bby_t, \xb_{t+1}) 
	+ 100 (1+\mu\gamma)\sum_{t=0}^{T-2}\beta_{t+1} L^2\gamma^2 \norm{\bbx_{t+1} - \bby_{t+1}}^2
	\end{aligned}
	\end{equation}

Next, we reformulate the inequality in Lemma~\ref{lem:main} as follows,
	\begin{align*}
		&\frac{\gamma}{\tau} \big(F(\bby_{t+1}) - F(x^*)\big) + \frac{1+\mu\gamma}{2m} \norm{\hzb_{t+1} - \mathbf{1}x^*}^2\\
		\leq&
		\frac{\gamma(1-\tau)}{\tau} \big( F(\bby_t) - F(x^*) \big) + \frac{1}{2m} \norm{\zb_t - \mathbf{1}x^* }^2 + \left(\frac{L\gamma}{2\tau} - \frac{1}{2\tau^2}\right) \norm{\bbx_{t+1} - \bby_{t+1}}^2 \\
		&- \frac{\gamma(1-\tau)}{\tau} D_f(\bby_t, \xb_{t+1})
		+ \frac{L\gamma}{2m\tau} \norm{\Pi \xb_{t+1}}^2 + \frac{1}{2m}\left( \gamma^2\norm{\Pi\bs_{t+1}}^2 + \norm{\Pi\hzb_{t+1}}^2  \right)\\
		\leq&
		\frac{\gamma(1-\tau)}{\tau} \big( F(\bby_t) - F(x^*) \big) + \frac{1}{2m} \norm{\hzb_t - \mathbf{1}x^* }^2 + \left(\frac{L\gamma}{2\tau} - \frac{1}{2\tau^2}\right) \norm{\bbx_{t+1} - \bby_{t+1}}^2 \\
		&- \frac{\gamma(1-\tau)}{\tau} D_f(\bby_t, \xb_{t+1})
		+ \frac{L\gamma}{2m\tau} \norm{\Pi \xb_{t+1}}^2 + \frac{1}{2m}\left( \gamma^2\norm{\Pi\bs_{t+1}}^2 + \norm{\Pi\hzb_{t+1}}^2  \right),
	\end{align*}
	where the last inequality is because of Lemma~\ref{lem:zz} with $\norm{\Pi\zb_t}^2 \leq \frac{1}{64} \norm{\Pi\hzb_t}^2$.

	By the choice of $\gamma$ and $\tau$, we have
	\begin{align*}
		\frac{\gamma}{\tau} \left(\frac{\gamma(1-\tau)}{\tau}\right)^{-1} = \frac{1}{1 - \tau} \ge 1 + \tau = 1 + \mu\gamma.
	\end{align*}
	Therefore, we can extract a factor of $(1+\mu\gamma)$ from the left.
	\begin{align*}
		&(1+\mu\gamma) \left( \frac{\gamma(1-\tau)}{\tau} \big(  F(\bby_{t+1}) - F(x^*)\big) + \frac{1}{2m} \norm{\hzb_{t+1} - \mathbf{1}x^*}^2 \right)\\
		\leq&
		\frac{\gamma(1-\tau)}{\tau} \big( F(\bby_t) - F(x^*) \big) + \frac{1}{2m} \norm{\hzb_t - \mathbf{1}x^* }^2 + \left(\frac{L\gamma}{2\tau} - \frac{1}{2\tau^2}\right) \norm{\bbx_{t+1} - \bby_{t+1}}^2 \\
		&- \frac{\gamma(1-\tau)}{\tau} D_f(\bby_t, \xb_{t+1})
	+ \frac{L\gamma}{2m\tau} \norm{\Pi \xb_{t+1}}^2 + \frac{1}{2m}\left( \gamma^2\norm{\Pi\bs_{t+1}}^2 + \norm{\Pi\hzb_{t+1}}^2  \right).
	\end{align*}
	Multiplying this inequality by $(1+\mu\gamma)^t$, summing over iterations $t=1$ to $t = T-1$, we can obtain the bound
	\begin{align*}
		&(1+\mu\gamma)^T \left( \frac{\gamma(1 -\tau)}{\tau} \big(F(\bby_T) - F(x^*)\big) + \frac{1}{2m} \norm{\hzb_T - \mathbf{1}x^*}^2 \right)\\
		\leq&
		\frac{\gamma(1-\tau)}{\tau} \big(F(\bby_1) - F(x^*)\big) + \frac{1}{2m} \norm{\hzb_1 - \mathbf{1}x^*}^2 \\
		&+ \sum_{t=1}^{T-1} (1+\mu\gamma)^t \left( \left(\frac{L\gamma}{2\tau} - \frac{1}{2\tau^2}\right) \norm{\bbx_{t+1} - \bby_{t+1}}^2  - \frac{\gamma(1-\tau)}{\tau} D_f(\bby_t, \xb_{t+1})  \right)\\
		&+ \sum_{t=1}^{T-1} (1+\mu\gamma)^t \left( \frac{L\gamma}{2m\tau} \norm{\Pi \xb_{t+1}}^2 + \frac{1}{2m}\left( \gamma^2\norm{\Pi\bs_{t+1}}^2 + \norm{\Pi\hzb_{t+1}}^2  \right)\right)\\
		\stackrel{\eqref{eq:rec}}{\leq}&
		\frac{\gamma(1-\tau)}{\tau} \big(F(\bby_1) - F(x^*)\big) + \frac{1}{2m} \norm{\hzb_1 - \mathbf{1}x^*}^2 \\
		&+ \sum_{t=0}^{T-1} (1+\mu\gamma)^t \left( \left(\frac{L\gamma}{2\tau} - \frac{1}{2\tau^2}\right) \norm{\bbx_{t+1} - \bby_{t+1}}^2 - \frac{\gamma(1-\tau)}{\tau} D_f(\bby_t, \xb_{t+1})  \right)\\
		&+\frac{2}{m\tau^2}  \left( \norm{\Pi\xb_1}^2 + c_1 \norm{\Pi\yb_1}^2 + c_2 \norm{\Pi \hzb_1}^2 + c_3\gamma^2 \norm{\Pi\bs_1}^2 \right) \cdot \sum_{t=0}^{T-1} (1 +\mu\gamma)^t \cdot 4^{-t} \\
		&+ 100 \sum_{t=1}^{T-1} (1+\mu\gamma)^t \cdot L\gamma^2 D_f(\bby_t, \xb_{t+1})  
		+ 100(1+\mu\gamma) \sum_{t=0}^{T-2} (1+\mu\gamma)^t L^2\gamma^2\norm{\bbx_{t+1} - \bby_{t+1}}^2 \\
		=&
		\frac{\gamma(1-\tau)}{\tau} \big(F(\bby_1) - F(x^*)\big) + \frac{1}{2m} \norm{\hzb_1 - \mathbf{1}x^*}^2 
		+ \sum_{t=1}^{T-1} (1+\mu\gamma)^t \cdot  \left(100L\gamma^2 - \frac{\gamma(1-\tau)}{\tau}\right) \cdot  D_f(\bby_t, \xb_{t+1}) \\
		&+ \sum_{t=1}^{T-2} \left( 100(1+\mu\gamma)L^2\gamma^2 + \frac{L\gamma}{2\tau} - \frac{1}{2\tau^2}\right)L^2\gamma^2\norm{\bbx_{t+1} - \bbx_{t+1}}^2
		\\
		& + 100 (1+\mu\gamma) \norm{\bbx_1 - \bby_1}^2 + (1+\mu\gamma)^T\left(\frac{L\gamma}{2\tau} - \frac{1}{2\tau^2}\right) \norm{\bbx_T - \bby_T}^2
		\\
		&+\frac{2}{m\tau^2}  \left( \norm{\Pi\xb_1}^2 + c_1 \norm{\Pi\yb_1}^2 + c_2 \norm{\Pi \hzb_1}^2 + c_3\gamma^2 \norm{\Pi\bs_1}^2 \right) \cdot \sum_{t=0}^{T-1} (1 +\mu\gamma)^t \cdot 4^{-t} \\
		\stackrel{\eqref{eq:tmg}}{\leq}&
		\frac{\gamma(1-\tau)}{\tau} \big(F(\bby_1) - F(x^*)\big) + \frac{1}{2m} \norm{\hzb_1 - \mathbf{1}x^*}^2 + (1+\mu\gamma) \norm{\bbx_1 - \bby_1}^2 \\
		&+\frac{2}{m\tau^2} \cdot\frac{1}{1 - \mu\gamma} \left( \norm{\Pi\xb_1}^2 + c_1 \norm{\Pi\yb_1}^2 + c_2 \norm{\Pi \hzb_1}^2 + c_3\gamma^2 \norm{\Pi\bs_1}^2 \right),
	\end{align*}
	where the last inequality is also because of the properly chosen $\gamma = \frac{1}{20\sqrt{\mu L}}$ and $\tau = \mu\gamma$ which implies that
	\begin{align*}
		&100L^2(1+\mu\gamma)\gamma^2 + \frac{L\gamma}{2\tau} - \frac{1}{2\tau^2} \leq \left(\frac{200}{400} + \frac{1}{2} - \frac{400}{2}\right) \frac{L}{\mu} < 0,\\
		&100L\gamma^2 - \frac{\gamma(1-\tau)}{\tau} = \frac{100}{400} \cdot \frac{1}{\mu} - \frac{1}{\mu} + \frac{1}{20\sqrt{\mu L}}  \stackrel{\mu\leq L}{\leq} \left( \frac{100}{400} + \frac{1}{20} - 1 \right) \frac{1}{\mu}<0.
	\end{align*}
	Thus, we can obtain that
	\begin{align*}
		\frac{1}{2m}\norm{\hzb_T -\mathbf{1}x^*}^2 
		\leq& 
		\frac{\gamma(1 -\tau)}{\tau} \big(F(\bby_T) - F(x^*)\big) + \frac{1}{2m} \norm{\hzb_T - \mathbf{1}x^*}^2 \\
		\leq&
		\left(1 + \mu\gamma\right)^{-T}  \cdot \Big( \frac{\gamma(1-\tau)}{\tau} \big(F(\bby_1) - F(x^*)\big) + \frac{1}{2m} \norm{\hzb_1 - \mathbf{1}x^*}^2 \\
		&+ \frac{2}{m\tau^2} \cdot \frac{1 }{1 - \mu\gamma} \cdot \left( \norm{\Pi\xb_1}^2 + c_1 \norm{\Pi\yb_1}^2 + c_2 \norm{\Pi \hzb_1}^2 + c_3\gamma^2 \norm{\Pi\bs_1}^2 \right) \Big)\\
		\leq&
		\left(1 + \frac{1}{20}\cdot \sqrt{\frac{\mu}{L}}\right)^{-T}  \cdot \Bigg( \frac{1}{\mu} \big(F(\bby_1) - F(x^*)\big) + \frac{1}{2m} \norm{\hzb_1 - \mathbf{1}x^*}^2 \\
		&+ \frac{ 3\cdot20^2 L}{m\mu}\cdot \left( \norm{\Pi\xb_1}^2 + c_1 \norm{\Pi\yb_1}^2 + c_2 \norm{\Pi \hzb_1}^2 + c_3\gamma^2 \norm{\Pi\bs_1}^2 \right)  \Bigg).
	\end{align*}
	Furthermore, by the fact that it holds that $(1 + a)^{-1} \leq 1 - \frac{a}{2}$ if $0\leq a \leq 1$ and multiply $2m$ both sides of above equation, we can obtain that
	\begin{align*}
		\norm{\hzb_T -\mathbf{1}x^*}^2
		\leq& 
		\left(1 - \frac{1}{40}\cdot \sqrt{\frac{\mu}{L}}\right)^{-T}  \cdot \Bigg( \frac{2m}{\mu} \big(F(\bby_1) - F(x^*)\big) +  \norm{\zb_1 - \mathbf{1}x^*}^2 \\
		&+ \frac{ 6\cdot20^2 L}{\mu}\cdot \left( \norm{\Pi\xb_1}^2 + c_1 \norm{\Pi\yb_1}^2 + c_2 \norm{\Pi \hzb_1}^2 + c_3\gamma^2 \norm{\Pi\bs_1}^2 \right)  \Bigg).
	\end{align*}
	Furthermore, since $\norm{\Pi\zb_T}^2 \leq \norm{\Pi\hzb_T}$, by Lemma~\ref{lem:zz}, we can obtain that $\norm{\zb_T-\mathbf{1}x^*}\leq \norm{\hzb_T-\mathbf{1}x^*}$. 
	It also holds that $\hzb_1 = \zb_1$.  
	Thus, we can obtain that
	\begin{align*}
	\norm{\zb_T -\mathbf{1}x^*}^2
	\leq& 
	\left(1 - \frac{1}{40}\cdot \sqrt{\frac{\mu}{L}}\right)^{-T}  \cdot \Bigg( \frac{2m}{\mu} \big(F(\bby_1) - F(x^*)\big) +  \norm{\zb_1 - \mathbf{1}x^*}^2 \\
	&+ \frac{ 2\cdot20^2 L}{\mu}\cdot \left( \norm{\Pi\xb_1}^2 + c_1 \norm{\Pi\yb_1}^2 + c_2 \norm{\Pi \zb_1}^2 + c_3\gamma^2 \norm{\Pi\bs_1}^2 \right)  \Bigg).
	\end{align*}
\end{proof}

\subsection{Proof of Theorems~\ref{thm:main1}}

\begin{lemma}
Suppose Assumption~\ref{ass:f} and Assumption~\ref{ass:g} hold with $\mu > 0$.	
The step size $\gamma_t$ and $\tau_t$ are set as Theorem~\ref{thm:main1}. 
Then it holds that
\begin{equation}\label{eq:aa}
\begin{aligned}
	&\sum_{t=1}^{T-1} \left[ \frac{L\gamma_t}{\tau_t}\norm{\Pi\xb_{t+1}}^2 + \gamma_t^2 \norm{\bs_{t+1}}^2 + \norm{\Pi\hzb_{t+1}}^2 \right]\\
	\leq&
	\frac{16\cdot 7}{9Lc_f}\left( \norm{\Pi\xb_1}^2 + \frac{9\tau_1^2}{32}\norm{\Pi\hzb_1}^2 + \frac{5}{14L^2c_f^2}\norm{\Pi\bs_1}^2 \right)\\
	&
	+ \frac{32}{9} \cdot   \sum_{t=1}^{T-1}\frac{\gamma_t}{\tau_t} \cdot\frac{50}{c_f} \cdot\left(  m D_f(\bby_t, \bbx_{t+1}) + 2 mL \norm{\bbx_{t+1} - \bby_{t+1}}^2  \right).
\end{aligned}
\end{equation}
\end{lemma}
\begin{proof}
	First, by the setting of the step size $\gamma_t = \frac{k+4}{2Lc_f}$ and the definition $\tau_t =\frac{1}{Lc_f\gamma_t} $, we have
	\begin{align*}
		L^2\gamma_t^2 <  \frac{L^2 \cdot 200^2 \dot \gamma_t^2}{100} =  \frac{L^2c_f^2 \gamma_t^2}{100} = \frac{1}{100\tau_t^2}.
	\end{align*} 
	Thus, the condition $L^2\gamma_t^2 \leq\frac{1}{100\cdot \tau_t^2}$ in Lemma~\ref{lem:err} holds.
	Furthermore, by the setting of $K$ in Theorem~\ref{thm:main}, the results in Lemma~\ref{lem:err} hold.
	
Next, using $\gamma_t \tau_t = \frac{1}{Lc_f}$,  we reformulate Eq.~\eqref{eq:xyzs} as follows:
\begin{equation}\label{eq:xyzs1}
\begin{aligned}
	&\norm{\Pi\xb_{t+1}}^2 + 64 \norm{\Pi \yb_{t+1}}^2 + \frac{9\tau_t^2}{32}\norm{\Pi\hzb_{t+1}}^2 + \frac{5\tau_t^2}{14}  \cdot \gamma_t^2   \norm{\Pi\bs_{t+1}}^2\\
	=&
	\norm{\Pi\xb_{t+1}}^2 + 64 \norm{\Pi \yb_{t+1}}^2 + \frac{9\tau_t^2}{32}\norm{\Pi\hzb_{t+1}}^2 + \frac{5}{14 L^2c_f^2}   \norm{\Pi\bs_{t+1}}^2
	\\
	\leq& 
	\frac{1}{4}\left( \norm{\Pi\xb_t}^2
	+
	64\norm{\Pi\yb_t}^2 +  \frac{9\tau_t^2}{32}\norm{\Pi\hzb_t}^2 +  \frac{5}{14L^2c_f^2}\norm{\Pi\hzb_t}^2 \right)\\
	&+\frac{25}{L^2c_f^2}\left(mL D_f (\bby_t, \xb_{t+1}) + mL^2 \norm{\bbx_t - \bby_t}^2\right)\\
	\leq&
	\frac{1}{4}\left( \norm{\Pi\xb_t}^2
	+
	64\norm{\Pi\yb_t}^2 +  \frac{9\tau_{t-1}^2}{32}\norm{\Pi\hzb_t}^2 +  \frac{5}{14L^2c_f^2}\norm{\Pi\hzb_t}^2 \right)\\
	&+\frac{25}{L^2c_f^2}\left(mL D_f (\bby_t, \xb_{t+1}) + mL^2 \norm{\bbx_t - \bby_t}^2\right),
\end{aligned}
\end{equation}
where the last inequality is because $\tau_k$ is non-increasing.

By denoting $\Delta_t = \frac{25}{L^2c_f^2}\left(mL D_f (\bby_t, \xb_{t+1}) + mL^2 \norm{\bbx_t - \bby_t}^2\right)$, we can obtain that
\begin{align*}
	&\sum_{t=1}^{T-1} \left[ \frac{L\gamma_t}{\tau_t}\norm{\Pi\xb_{t+1}}^2 + \gamma_t^2 \norm{\bs_{t+1}}^2 + \norm{\Pi\hzb_{t+1}}^2 \right] \\
	=&
	\sum_{t=1}^{T-1}\frac{L\gamma_t}{\tau_t} \left[ \norm{\Pi\xb_{t+1}}^2 +  \frac{\tau_t}{L\gamma_t} \frac{32}{9\tau_t^2} \cdot \frac{9\tau_t^2}{32}\norm{\Pi\hzb_{t+1}}^2 +  \frac{\tau_t}{L\gamma_t} \frac{14}{5\tau_t^2}  \cdot  \frac{5\tau_t^2}{14} \gamma_t^2 \norm{\Pi\bs_{t+1}}^2 +  \right]\\
	=&
	\sum_{t=1}^{T-1}\frac{L\gamma_t}{\tau_t} \left[ \norm{\Pi\xb_{t+1}}^2 +   \frac{32}{9L\gamma_t\tau_t} \cdot \frac{9\tau_t^2}{32}\norm{\Pi\hzb_{t+1}}^2 +   \frac{14}{5L\gamma_t\tau_t}  \cdot  \frac{5\tau_t^2}{14} \gamma_t^2 \norm{\Pi\bs_{t+1}}^2  \right]\\
	=&
	\sum_{t=1}^{T-1}\frac{L\gamma_t}{\tau_t} \left[ \norm{\Pi\xb_{t+1}}^2 +   \frac{32c_f}{9} \cdot \frac{9\tau_t^2}{32}\norm{\Pi\hzb_{t+1}}^2 +   \frac{14c_f}{5}  \cdot  \frac{5\tau_t^2}{14} \gamma_t^2 \norm{\Pi\bs_{t+1}}^2  \right]\\
	\leq&
	\frac{32c_f}{9} \cdot \sum_{t=1}^{T-1}\frac{L\gamma_t}{\tau_t} \left[ \norm{\Pi\xb_{t+1}}^2 +    \frac{9\tau_t^2}{32}\norm{\Pi\hzb_{t+1}}^2 +    \frac{5}{14L^2c_f^2}  \norm{\Pi\bs_{t+1}}^2  \right]\\
	\stackrel{\eqref{eq:xyzs1}}{\leq}&
	 \frac{32c_f}{9}\cdot \sum_{t=1}^{T-1}\frac{L\gamma_t}{\tau_t}\left(\frac{1}{4}\right)^t \left( \norm{\Pi\xb_1}^2 + \frac{9\tau_1^2}{32}\norm{\Pi\hzb_1}^2 + \frac{5}{14L^2c_f^2}\norm{\Pi\bs_1}^2 \right)+ \frac{32c_f}{9} \cdot\sum_{t=1}^{T-1} \frac{L\gamma_t}{\tau_t} \sum_{k=1}^{t}\left(\frac{1}{4}\right)^{t-k} \Delta_k.
\end{align*}

Furthermore,
\begin{align*}
	\sum_{t=1}^{T-1}\frac{\gamma_t}{\tau_t}\left(\frac{1}{4}\right)^t 
	= 
	\sum_{t=1}^{T-1}  Lc_f \gamma_t^2 \left(\frac{1}{4}\right)^{t}
	=
	\sum_{t=1}^{T-1} Lc_f \frac{(t+4)^2}{4L^2c_f^2}\left(\frac{1}{4}\right)^t
	\leq 
	\frac{7}{2L^2c_f^2},
\end{align*}
where the last inequality is because of Lemma~\ref{lem:tsum}.

We also have
\begin{equation}\label{eq:xx}
	\frac{\gamma_{t+1}}{\tau_{t+1}} \cdot \frac{1}{4} = \frac{(t+5)^2}{4Lc_f} \cdot \frac{1}{4} \leq \frac{(t+4)^2}{4Lc_f} \cdot \frac{1}{2} = \frac{\gamma_t}{\tau_t} \cdot \frac{1}{2}.
\end{equation}
Thus, we can obtain that 
\begin{align*}
	\sum_{t=1}^{T-1} \frac{L\gamma_t}{\tau_t} \sum_{k=1}^{t}\left(\frac{1}{4}\right)^{t-k} \Delta_k
	\stackrel{\eqref{eq:xx}}{\leq} 
	\sum_{t=1}^{T-1}  \sum_{k=1}^{t} \frac{L\gamma_k}{\tau_k}\left(\frac{1}{2}\right)^{t-k} \Delta_k
	\leq
	2L\sum_{t=1}^{T-1}\frac{\gamma_t}{\tau_t} \Delta_t,
\end{align*}
where the last inequality is because of Lemma~\ref{lem:recs}.

Combining above results, we obtain
\begin{align*}
&\sum_{t=1}^{T-1} \left[ \frac{L\gamma_t}{\tau_t}\norm{\Pi\xb_{t+1}}^2 + \gamma_t^2 \norm{\bs_{t+1}}^2 + \norm{\Pi\hzb_{t+1}}^2 \right]\\
\leq&
\frac{16\cdot 7}{9Lc_f}\left( \norm{\Pi\xb_1}^2 + \frac{9\tau_1^2}{32}\norm{\Pi\hzb_1}^2 + \frac{5}{14L^2c_f^2}\norm{\Pi\bs_1}^2 \right)\\
&
+ \frac{32}{9} \cdot   \sum_{t=1}^{T-1}\frac{\gamma_t}{\tau_t} \cdot\frac{50}{c_f} \cdot\left(  m D_f(\bby_t, \bbx_{t+1}) +  mL \norm{\bbx_{t+1} - \bby_{t+1}}^2  \right)\\
\leq&
\frac{16\cdot 7}{9Lc_f}\left( \norm{\Pi\xb_1}^2 + \frac{9\tau_1^2}{32}\norm{\Pi\hzb_1}^2 + \frac{5}{14L^2c_f^2}\norm{\Pi\bs_1}^2 \right)\\
&
+ \frac{32}{9} \cdot   \sum_{t=1}^{T-1}\frac{\gamma_t}{\tau_t} \cdot\frac{50}{c_f} \cdot\left(  m D_f(\bby_t, \bbx_{t+1}) +  2mL \norm{\bbx_{t+1} - \bby_{t+1}}^2  \right),
\end{align*}
where the last inequality is because of
\begin{align*}
	&\sum_{t=1}^{T-1} \frac{\gamma_t}{\tau_t} \norm{\bbx_t - \bby_t}^2 
	= 
	\frac{\gamma_1}{\tau_1} \norm{\bbx_1 - \bby_1}^2 + \sum_{t=2}^{T-1}  \frac{\gamma_t}{\tau_t} \norm{\bbx_t - \bby_t}^2\\
	=&
	\sum_{t=1}^{T-2} \frac{\gamma_{t+1}}{\tau_{t+1}} \cdot \left(\frac{\gamma_t}{\tau_t}\right)^{-1} \cdot \frac{\gamma_t}{\tau_t} \norm{\bbx_{t+1} - \bby_{t+1}}^2
	\leq
	2 \sum_{t=1}^{T-1}\frac{\gamma_t}{\tau_t} \norm{\bbx_{t+1} - \bby_{t+1}}^2.
\end{align*}
\end{proof}

Combining with results in previous sections, we can prove Theorem~\ref{thm:main1}.

\begin{proof}[Proof of Theorem~\ref{thm:main1}]
By setting $\mu = 0$ in Lemma~\ref{lem:main} and by Lemma~\ref{lem:zz}, we can obtain 
\begin{align*}
0 \leq& \frac{\gamma_t(1-\tau_t)}{\tau_t} F(\bby_t) - \frac{\gamma_t}{\tau_t}F(\bby_{t+1}) + \gamma_t F(x^*) + \left(\frac{L\gamma_t}{2\tau_t} - \frac{1}{2\tau_t^2}\right) \norm{\bbx_{t+1} - \bby_{t+1}}^2 \\
&+\frac{1}{2m} \norm{\hzb_t - \mathbf{1}x^*}^2 - \frac{1}{2m} \norm{\hzb_{t+1} - \mathbf{1}x^*}^2 - \frac{\gamma_t(1-\tau_t)}{\tau_t} D_f(\bby_t, \xb_{t+1})\\
&+ \frac{L\gamma_t}{2m\tau_t} \norm{\Pi \xb_{t+1}}^2 + \frac{1}{2m}\left( \gamma_t^2\norm{\Pi\bs_{t+1}}^2 + \norm{\Pi\hzb_{t+1}}^2 \right).
\end{align*}
Summing above result over the iteration $k=1$ to $k = T-1$, we can obtain that
\begin{align*}
	0\leq& \frac{1}{2m}\norm{\hzb_1 - \mathbf{1}x^*}^2 - \frac{1}{2m}\norm{\hzb_T - \mathbf{1}x^*}^2 + \sum_{t=1}^{T-1} \left[ \frac{\gamma_t(1-\tau_t)}{\tau_t} F(\bby_t) - \frac{\gamma_t}{\tau_t} F(\bby_{t+1}) + \gamma_t F(x^*) \right]\\
	&+\sum_{t=1}^{T-1} \left[ \frac{\gamma_t}{\tau_t}\left(\frac{L}{2} - \frac{1}{2\tau_t\gamma_t}\right) \norm{\bbx_{t+1} - \bby_{t+1}}^2 - \frac{\gamma_t(1-\tau_t)}{\tau_t} D_f(\bby_t, \xb_{t+1})\right]\\
	&+\frac{1}{2m}\sum_{t=1}^{T-1} \left[ \frac{L\gamma_t}{\tau_t}\norm{\Pi\xb_{t+1}}^2 + \gamma_t^2 \norm{\bs_{t+1}}^2 + \norm{\Pi\hzb_{t+1}}^2 \right].
\end{align*}

Combining with Eq.~\eqref{eq:aa}, we can obtain that
\begin{align*}
	0\leq& 
	\frac{1}{2m}\norm{\hzb_1 - \mathbf{1}x^*}^2 - \frac{1}{2m}\norm{\hzb_T - \mathbf{1}x^*}^2 + \sum_{t=1}^{T-1} \left[ \frac{\gamma_t(1-\tau_t)}{\tau_t} F(\bby_t) - \frac{\gamma_t}{\tau_t} F(\bby_{t+1}) + \gamma_t F(x^*) \right]\\
	&+\sum_{t=1}^{T-1} \left[ \frac{\gamma_t}{\tau_t}\left( \frac{32}{9} \cdot \frac{100}{c_f} L +\frac{L}{2} - \frac{1}{2\tau_t\gamma_t}\right) \norm{\bbx_{t+1} - \bby_{t+1}}^2 + \frac{\gamma_t}{\tau_t} \left( \frac{32}{9} \cdot \frac{25}{c_f} -1+\tau_t \right) D_f(\bby_t, \xb_{t+1})\right]\\
	&
	+\frac{16\times 7}{9Lc_f} \left(\norm{\Pi\xb_1}^2 + \frac{9\tau_1^2}{32}\norm{\Pi\hzb_1}^2 + \frac{5}{14L^2c_f^2}\norm{\Pi\bs_1}^2 \right)\\
	\leq&
	\frac{1}{2m}\norm{\hzb_1 - \mathbf{1}x^*}^2  + \sum_{t=1}^{T-1} \left[ \frac{\gamma_t(1-\tau_t)}{\tau_t} F(\bby_t) - \frac{\gamma_t}{\tau_t} F(\bby_{t+1}) + \gamma_t F(x^*) \right]\\
	&
	+\frac{16\times 7}{9Lc_f} \left(\norm{\Pi\xb_1}^2 + \frac{9\tau_1^2}{32}\norm{\Pi\hzb_1}^2 + \frac{5}{14L^2c_f^2}\norm{\Pi\bs_1}^2 \right)\\
	=&
	\frac{1}{2m}\norm{\hzb_1 - \mathbf{1}x^*}^2  + \sum_{t=1}^{T-1} \Big[ \left(Lc_f \gamma_t^2 - \gamma_t\right) F(\bby_t) - Lc_f\gamma_t^2 \cdot F(\bby_{t+1}) + \gamma_t F(x^*) \Big]\\
	&
	+\frac{16\times 7}{9Lc_f} \left(\norm{\Pi\xb_1}^2 + \frac{9\tau_1^2}{32}\norm{\Pi\hzb_1}^2 + \frac{5}{14L^2c_f^2}\norm{\Pi\bs_1}^2 \right),
\end{align*}
where the second inequality is because of  the following facts
\begin{align*}
	&\frac{32}{9} \cdot \frac{100}{c_f} L +\frac{L}{2} - \frac{1}{2\tau_t\gamma_t} 
	=
	L\left(\frac{3200}{9c_f} + \frac{1}{2} - \frac{c_f}{2} \right)
	=
	L\left(\frac{16}{9} + \frac{1}{2} - 100\right)<0,\\
	&\frac{32}{9} \cdot \frac{25}{c_f} - 1 + \tau_t 
	= 
	\frac{32\cdot 25}{9c_f} - 1 + \frac{2}{t+4}
	\leq
	\frac{32\cdot 25}{9c_f} - 1+ \frac{2}{5} = \frac{32\cdot 25}{9\cdot 200} - \frac{3}{5}<0.
\end{align*}

Furthermore, it holds that
\begin{align*}
	&Lc_f \gamma_t^2 - \gamma_t - \left(Lc_f\gamma_{t-1}^2 - \frac{1}{4Lc_f}\right)
	=
	Lc_f (\gamma_t - \gamma_{t-1}) (\gamma_t+\gamma_{t-1}) - \gamma_t + \frac{1}{4Lc_f}\\
	=&
	\frac{\gamma_t+\gamma_{t-1}}{2} - \gamma_t + \frac{1}{4Lc_f}
	=
	\frac{\gamma_{t-1} - \gamma_t}{2} + \frac{1}{4Lc_f} = 0.
\end{align*}
Thus, we can conclude that
\begin{align*}
	Lc_f \gamma_t^2 - \gamma_t = Lc_f\gamma_{t-1}^2 - \frac{1}{4Lc_f}.
\end{align*}
Thus,
\begin{align*}
&\sum_{t=1}^{T-1} \Big[ \left(Lc_f \gamma_t^2 - \gamma_t\right) F(\bby_t) - Lc_f\gamma_t^2 \cdot F(\bby_{t+1}) + \gamma_t F(x^*) \Big]\\
=&
\sum_{t=1}^{T-1} \Big[ \left(Lc_f \gamma_{t-1}^2 - \frac{1}{4Lc_f}\right) F(\bby_t) - Lc_f\gamma_t^2 \cdot F(\bby_{t+1}) + \Big(Lc_f\gamma_t^2 - Lc_f\gamma_{t-1}^2 + \frac{1}{4Lc_f}\Big) F(x^*) \Big]\\
=&
-Lc_f\gamma_{T-1}^2 \Big( F(\bby_T) - F(x^*) \Big)  - \frac{1}{4Lc_f}\sum_{t=1}^{T-1} \Big(F(\bby_t) - F(x^*)\Big)
+ \big(Lc_f\gamma_1^2- \gamma_1\big) \Big(F(\bby_1) - F(x^*)\Big).
\end{align*}

Using the facts that $Lc_f\gamma_{T-1}^2 = \frac{(T+3)^2}{4Lc_f}$, $Lc_f\gamma_1^2 - \gamma_1 = \frac{15}{4Lc_f}$
, and $F(\bby_t) \leq
F(x^*)$, we have
\begin{align*}
0 \leq&
\frac{1}{2m}\norm{\hzb_1 - \mathbf{1}x^*}^2  + \sum_{t=1}^{T-1} \Big[ \left(Lc_f \gamma_t^2 - \gamma_t\right)  F(\bby_t) - Lc_f\gamma_t^2 \cdot F(\bby_{t+1}) + \gamma_t F(x^*) \Big]\\
&
+\frac{16\times 7}{9Lc_f} \left(\norm{\Pi\xb_1}^2 + \frac{9\tau_1^2}{32}\norm{\Pi\hzb_1}^2 + \frac{5}{14L^2c_f^2}\norm{\Pi\bs_1}^2 \right)\\
\leq&
\frac{1}{2m}\norm{\hzb_1 - \mathbf{1}x^*}^2 -Lc_f\gamma_{T-1}^2 \Big( F(\bby_T) - F(x^*) \Big)- \frac{1}{4Lc_f}\sum_{t=1}^{T-1} \Big(F(\bby_t) - F(x^*)\Big)
\\
&+ \big(Lc_f\gamma_1^2- \gamma_1\big) \Big(F(\bby_1) - F(x^*)\Big)
+\frac{16\times 7}{9Lc_f} \left(\norm{\Pi\xb_1}^2 + \frac{9\tau_1^2}{32}\norm{\Pi\hzb_1}^2 + \frac{5}{14L^2c_f^2}\norm{\Pi\bs_1}^2 \right)\\
\leq&
\frac{1}{2m}\norm{\hzb_1 - \mathbf{1}x^*}^2 - \frac{(T+3)^2}{4Lc_f} \Big( F(\bby_T) - F(x^*) \Big) + \frac{15}{4Lc_f}\Big(F(\bby_1) - F(x^*)\Big) \\
&
+ \frac{16\cdot 7}{9Lc_f} \left(\norm{\Pi\xb_1}^2 + \frac{9\tau_1^2}{32}\norm{\Pi\hzb_1}^2 + \frac{5}{14L^2c_f^2}\norm{\Pi\bs_1}^2 \right).
\end{align*}
Above equation implies that
\begin{align*}
	\frac{(T+3)^2}{4Lc_f} \Big( F(\bby_T) - F(x^*) \Big) 
	\leq& 
	\frac{1}{2m}\norm{\hzb_1 - \mathbf{1}x^*}^2 + \frac{15}{4Lc_f} \Big(F(\bby_1) - F(x^*)\Big) \\
	&
	+ \frac{16\cdot 7}{9Lc_f} \left(\norm{\Pi\xb_1}^2 + \frac{9\tau_1^2}{32}\norm{\Pi\hzb_1}^2 + \frac{5}{14L^2c_f^2}\norm{\Pi\bs_1}^2 \right).
\end{align*}
Dividing $\frac{(T+3)^2}{4Lc_f}$ for both sides of above equation, we can obtain
\begin{align*}
	F(\bby_T) - F(x^*)  
	\leq& 
	\frac{15}{(T+3)^2} \Big(F(\bby_1) - F(x^*)\Big) + \frac{2Lc_f}{m(T+3)^2} \norm{\hzb_1 - \mathbf{1}x^*}^2 \\
	&
	+ \frac{50}{(T+3)^2} \left(\norm{\Pi\xb_1}^2 + \frac{9\tau_1^2}{32}\norm{\Pi\hzb_1}^2 + \frac{5}{14L^2c_f^2}\norm{\Pi\bs_1}^2 \right).
\end{align*}
Combining with the fact $\hzb_1 = \zb_1$, we can conclude the proof. 
\end{proof}
\section{Experiments}

In the previous sections, we have given the theoretical analysis of {our algorithm}. 
In this section, we will validate the effectiveness and computational efficiency of our algorithm empirically. 
We will conduct experiments on the sparse logistic regression problem whose objective function satisfies the form~\eqref{eq:prob} with 
\begin{equation}
		f_i(x) = \frac{1}{n}\sum_{j=1}^{n} \log \big(1+\exp(-b_{ij}\langle a_{ij}, x\rangle)\big), \qquad\mbox{and}\qquad
		g(x) = \sigma\norm{x}_1 + \frac{\mu}{2}\norm{x}^2, \label{eq:lgr}
\end{equation}
where $a_{ij} \in \RR^{d}$ and $b_{ij}\in\{-1,1\}$ are the $j$-th input vector and the corresponding label on the $i$-th agent. 
We can observer that each $f_i(x)$ in Eq.~\eqref{eq:lgr} is convex and smooth. 
At the same time, $g(x)$ is $\mu$-strongly convex.  
Thus, the sparse logistic regression problem satisfies the assumptions required in our algorithm.

\begin{figure*}[!ht]
	\begin{center}
		\centering
		\subfigure[\textsf{a9a error versus gradient computation.}]{\includegraphics[width=60mm]{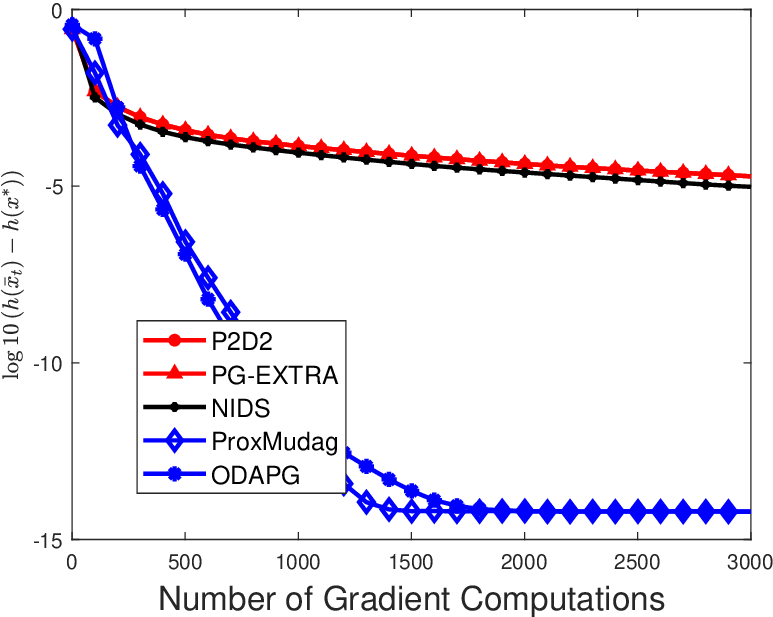}}~
		\subfigure[\textsf{a9a error versus communication.}]{\includegraphics[width=60mm]{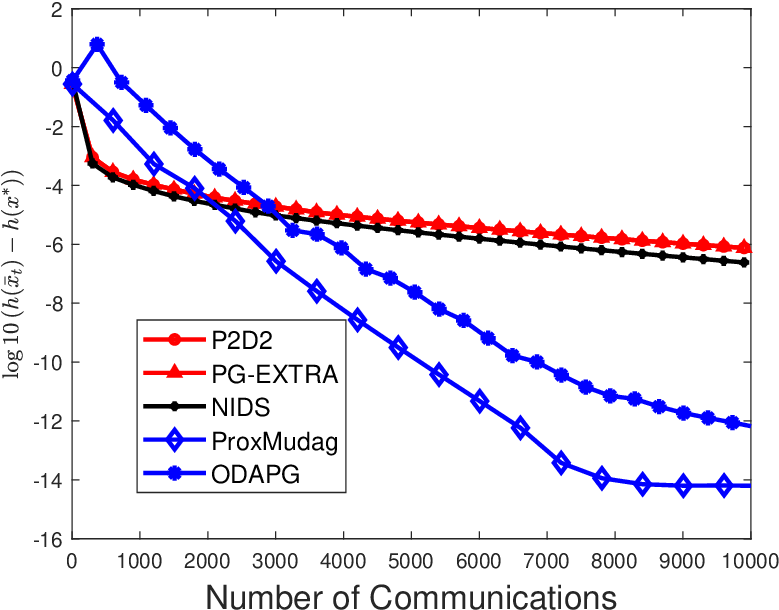}}\\
		\subfigure[\textsf{w8a error versus gradient computation.}]{\includegraphics[width=60mm]{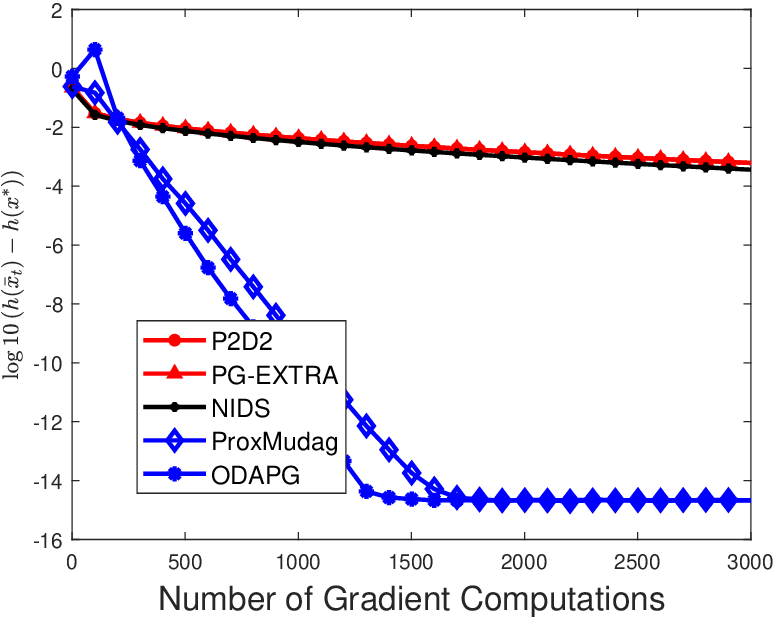}}~
		\subfigure[\textsf{w8a error versus communication.}]{\includegraphics[width=60mm]{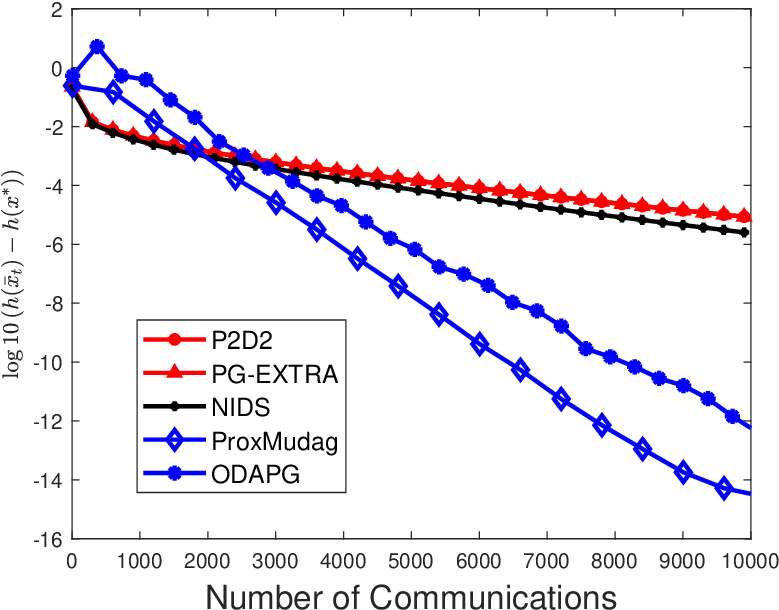}}
	\end{center}
	\caption{Algorithm evaluation with comparison to the decentralized optimization methods on the network with $1 - \lambda_2(W) = 0.05$.}
	\label{fig:result}
\end{figure*}

\paragraph{Experiments Setting}
In our experiments, we consider random networks where each pair of agents have a connection with a probability of $p = 0.1$. 
We set $W = I - \frac{\mathbf{L}}{\lambda_1(\mathbf{L})}$ where $\mathbf{L}$ is the Laplacian matrix associated with a weighted graph,
and $\lambda_1(\mathbf{L})$ is the largest eigenvalue of $\mathbf{L}$. 
We set $m = 100$, that is, there exists $100$
agents in this network. 
In our experiments, the gossip matrix $W$ satisfies  $1-\lambda_2(W) = 0.05$.

We conduct experiments on the datasets `w8a' and `w9a' which can be downloaded in libsvm datasets.
For `w8a', we set $n = 497$ and $d = 300$.
For `a9a', we set $n = 325$ and $d = 123$.
We set $\sigma = 10^{-4}$ for all datasets  and set $\mu$ as  $10^{-4}$ which leads to a large condition number of the objective function.

\paragraph{Comparison with Existing Works}
We compare our work with state-of-the-art algorithms \texttt{PG-EXTRA} \citep{ShiLWY15}, \texttt{NIDS} \citep{li2019decentralized}, Decentralized Proximal Algorithm (\texttt{DPA}) \citep{AlghunaimYS19} and \texttt{ProxMudag} \citep{Ye2023}. 
Both our algorithm and \texttt{ProxMudag} need ``FastMix'' to achieve consensus. 
In our experiments, we set $K=3$ in the ``FastMix'' for our algorithm and \texttt{ProxMudag}.
In our experiments, the parameters of all algorithms are well-tuned. 
We report experiment results in Figure~\ref{fig:result}.
We can observe that \texttt{ODAPG} and \texttt{ProxMudag} take much less computational cost than other algorithms because these two algorithms use Nesterov's acceleration to achieve faster convergence rates. That is,
the computation complexity of \texttt{ODAPG} and \texttt{ProxMudag} is linear to $\sqrt{L/\mu}$ while the computation complexities of other algorithms without the acceleration such as \texttt{NIDS} are linear to $L/\mu$.
This matches our theoretical analysis of \texttt{ODAPG} over the computation complexity.
\texttt{ODAPG} also shows great advantages over other  decentralized proximal algorithms without the acceleration on the communication cost.
Though \texttt{ODAPG} takes three times of local communication while other algorithms communicate only once for each iteration, 
\texttt{ODAPG} still requires much less communication costs because of its fast convergence rate.
Compared with \texttt{ProxMudag}, \texttt{ODAPG} achieves a similar computation cost with the one of \texttt{ProxMudag}. 
This is because both these two algorithms apply Nesterov's acceleration to promote the convergence rates.
For the communication cost, our \texttt{ODAPG} takes a little more communication cost to achieve the same precision solution. 
This is because \texttt{ODAPG} takes three ``FastMix'' steps while \texttt{ProxMudag} takes two ``FastMix'' steps for each iteration of algorithms.
However, theoretically,  \texttt{ODAPG} only requires $K = \cO(1)$ for the  ``FastMix'' step while \texttt{ProxMudag} requires $K = \cO\left(\log \frac{L}{\mu}\right)$. 

\section{Conclusion and Future Work}

In this paper, we have proposed the  first optimal algorithm for the decentralized convex composite optimization problem.
When $f_i(x)$ is $L$-smooth and $g(x)$ is $\mu$-strong convex, our algorithm can achieve a computation complexity $\cO\left(\sqrt{\frac{L}{\mu}}\log\frac{1}{\epsilon}\right)$ and a communication complexity $\cO\left(\sqrt{\frac{L}{\mu(1 - \lambda_2(W))}}\log\frac{1}{\epsilon}\right)$. 
When $f(x)$ is $L$-smooth and $g(x)$ is generally convex, our algorithm can achieve a computation complexity $\cO\left(\sqrt{\frac{1}{\epsilon}}\right)$ and a communication complexity $\cO\left(\sqrt{\frac{1}{\epsilon(1 - \lambda_2(W))}}\right)$.
That is, our algorithm achieves both the optimal computation and communication  complexities when  $f_i(x)$ is smooth and the regularization term $g(x)$ is (strongly)-convex. 
For the problem that that $f_i(x)$ is $L$-smooth and $\mu$-strongly convex while $g(x)$ is convex, our algorithm can be easily extended to solve this kind of problem and also achieve the optimal computation and communication complexities. 
Our experiments also validate the effectiveness and efficiency of our algorithm both in  computation and communication.

Our algorithm requires the multi-consensus which takes an inner communication loop to reduce the consensus error. \citet{song2023optimal} propose the optimal decentralized algorithm for smooth functions with only single-loop communication.
We conjecture our algorithm can be implemented similarly by the technique in \citet{song2023optimal} to remove the multi-consensus steps in our algorithm.

\appendix

\section{Some Useful Lemmas}

\begin{lemma}
	Letting $x$ and $y$ be two $d$-dimensional vectors, then it holds that
	\begin{equation}
		\dotprod{x, y} \leq \norm{x} \norm{y} \leq \frac{\norm{x}^2 + \norm{y}^2}{2}, \label{eq:cauchy}
	\end{equation}
	and
	\begin{equation}
		\norm{x + y}^2 \leq 2 \left(\norm{x}^2 + \norm{y}^2\right). \label{eq:xy}
	\end{equation}
\end{lemma}
\begin{proof}
	First, by the Cauchy's inequality, we have $\dotprod{x,y} \leq \norm{x}\norm{y}$. Using the fact $2ab \leq a^2 + b^2$ for any scalars $a$ and $b$, we can obtain that $\norm{x}\norm{y} \leq \frac{\norm{x}^2 + \norm{y}^2}{2}$.
	
	Furthermore, 
	\begin{align*}
		\norm{x + y}^2 = \norm{x}^2 + \norm{y}^2 + 2 \dotprod{x,y} \stackrel{\eqref{eq:cauchy}}{\leq} 2 \left(\norm{x}^2 + \norm{y}^2\right).
	\end{align*}
\end{proof}

\begin{lemma}[Lemma 7 of \citet{driggs2021accelerating}]\label{lem:recs}
	Given a non-negative sequence $\sigma_t$, a constant $\rho \in (0, 1]$, and an index $T \geq 1$, the following estimate holds:
	\begin{equation*}
		\sum_{t=1}^{T} \sum_{\ell=1}^{t} (1-\rho)^{t - \ell} \sigma_\ell \leq \frac{1}{\rho}\sum_{t=1}^{T}\sigma_t.
	\end{equation*}
\end{lemma}

\begin{lemma}
If $0<\mu\gamma <1$, given $T > 1$, then it holds that
\begin{equation}
\sum_{t=0}^{T-1} (1 + \mu\gamma)^t \cdot 4^{-t} \leq \frac{1}{1 - \mu\gamma}. \label{eq:tmg}
\end{equation}
\end{lemma}
\begin{proof}
	It holds that
\begin{align*}
\sum_{t=0}^{T-1} (1 + \mu\gamma)^t \cdot 4^{-t} = \frac{1 - (1+\mu\gamma)^T \cdot 4^{-T}}{1 - \frac{1+\mu\gamma}{4}} \leq \frac{1}{1 - \mu\gamma}.
\end{align*}
\end{proof}

\begin{lemma}\label{lem:tsum}
	It holds that $\sum_{t=1}^{T-1} (t+4)^2 \left(\frac{1}{4}\right)^t \leq 14$.
\end{lemma}
\begin{proof}
	Letting us denote $g(t) =(t+4)^2\left(\frac{1}{4}\right)^t$, then its derivative $g'(t)$ satisfies 
	\begin{align*}
		 g'(t) = 2(t+4)\left(\frac{1}{4}\right)^t - (t+4)^2\left(\frac{1}{4}\right)^t\ln(4) < 0.
	\end{align*}
	This implies that $g(t)$ is a deceasing function.
	Thus, we have
	\begin{align*}
		\sum_{t=1}^{T-1} (t+4)^2 \left(\frac{1}{4}\right)^t 
		\leq&
		\int_{1}^{T}(t+4)^2\left(\frac{1}{4}\right)^t\;dt 
		\leq 
		\int_{1}^{T}(t+4)^2 e^{-t} \;dt\\
		=&
		-e^{-T}(T+4)^2 + 25/e + 2\int_{1}^{T} e^{-t} (t+4)\;dt\\
		=&
		-e^{-T}(T+4)^2 + 25/e + 2 \left(-e^{-T}(T+4) + 5/e + \int_{1}^{T} e^{-t}\;dt\right)\\
		\leq&
		\frac{37}{e} \leq 14.
	\end{align*}
	
\end{proof}
\begin{lemma}[Lemma 3 of \cite{driggs2021accelerating}]\label{lem:g}
	Suppose $g$ is $\mu$-strongly convex with $\mu \ge 0$, and suppose $z = \prox_{\gamma g} (x - \gamma d)$ for some $x, d \in\RR^d$ and constant $\gamma$. Then for $y \in\RR^d$,
	\begin{equation}
		\gamma\dotprod{d, z - y} 
		\leq
		\frac{1}{2}\norm{x - y}^2 - \frac{1+\mu\gamma}{2} \norm{z - y}^2 - \frac{1}{2}\norm{z - x} - \gamma g(z) + \gamma g(y). 
	\end{equation}
\end{lemma}

\begin{lemma}
	If $f_i(\cdot)$ is convex and $L$-smooth, then it holds that for $x, y\in \RR^d$
	\begin{align}
		\norm{\nabla f_i(x) - \nabla f_i(y)}^2 \le& 2L \big( f_i(x) - f_i(y) - \dotprod{\nabla f_i(y), x - y} \big), \label{eq:L_breg} \\
		\norm{\nabla f_i(x) - \nabla f_i(y)} \leq& L \norm{x - y}. \label{eq:L_n}
	\end{align} 
\end{lemma}

\begin{lemma}\label{lem:zz}
Letting $K$ in Eq.~\eqref{eq:hz} be chosen large enough that $\norm{\Pi\zb_t}^2 \leq \rho^2 \norm{\Pi\hzb_t}^2 $ with $\rho^2 \leq 1$, then it holds that
\begin{equation}
	\norm{\zb_t - \mathbf{1}x^*}^2 \leq \norm{\hzb_t - \mathbf{1}x^*}^2.
\end{equation}
\end{lemma}
\begin{proof}
	First, we have
\begin{align*}
	&\norm{\zb_t - \mathbf{1}x^*}^2
	= \norm{\zb_t - \hzb_t + \hzb_t - \mathbf{1}x^*}^2\\
	=& 
	\norm{\zb_t - \hzb_t}^2 + \norm{\hzb_t - \mathbf{1}x^*}^2
	+2\sum_{i=1}^{m}\dotprod{\zbti - \hzb_t^{(i)}, \hzb_t^{(i)} - x^*}\\
	=&
	\norm{\zb_t - \hzb_t}^2 + \norm{\hzb_t - \mathbf{1}x^*}^2 + 2\sum_{i=1}^{m}\dotprod{\zbti - \hzb_t^{(i)}, \hzb_t^{(i)}}\\
	=&
	\norm{\zb_t}^2 + \norm{\hzb_t}^2 -  2\sum_{i=1}^{m}\dotprod{\zbti, \hzb_t^{(i)}}  + \norm{\hzb_t - \mathbf{1}x^*}^2 + 2\sum_{i=1}^{m}\dotprod{\zbti - \hzb_t^{(i)}, \hzb_t^{(i)}}\\
	=&
	\norm{\zb_t}^2 - \norm{\hzb_t}^2 + \norm{\hzb_t - \mathbf{1}x^*}^2.
\end{align*}
Furthermore,
\begin{align*}
&\norm{\zb_t}^2 - \norm{\hzb_t}^2
=
\norm{\zb_t}^2 + \norm{\mathbf{1}\bbz_t}^2 - 2\sum_{i=1}^{m}\dotprod{\zbti, \bbz_t} -\left(\norm{\hzb_t}^2 + \norm{\mathbf{1}\bbz_t}^2 - 2\sum_{i=1}^{m}\dotprod{\hzb_t^{(i)}, \bbz_t}\right)\\
=&
\norm{\zb_t - \mathbf{1}\bbz_t}^2 - \norm{\hzb_t - \mathbf{1}\bbz_t}^2
\leq 0,
\end{align*}
where the last inequality is because of the assumption $\norm{\zb_t - \mathbf{1}\bbz_t}^2 \leq  \norm{\hzb_t - \mathbf{1}\bbz_t}^2$.
Combining above two equations, we can obtain the result.
\end{proof}

\bibliography{ref.bib}
\bibliographystyle{apalike2}

\end{document}